\documentclass[12pt]{article}
\usepackage[backref=page,colorlinks,bookmarksopen=true]{hyperref}
\usepackage{pdfsync}
\usepackage{amssymb}
\usepackage{latexsym}
\usepackage{amsmath}
\usepackage[all]{xy}
\newtheorem{theorem}{Theorem}[section]
\newtheorem{lemma}[theorem]{Lemma}
\newtheorem{proposition}[theorem]{Proposition}
\newtheorem{corollary}[theorem]{Corollary}
\newtheorem{definition}[theorem]{Definition}
\newtheorem{example}[theorem]{Example}

\newtheorem{remark}[theorem]{Remark}

\begin{document}
\title{Yetter-Drinfel'd algebras and coideals of Weak Hopf $C^*$-Algebras}
\author{Leonid Vainerman Jean-Michel Vallin}
\date{29 June 2020}
\maketitle
%127317
\footnote{LMNO Universit\'e de Caen, \ IDP Universit\'e d'Orl\'eans/IMJ Universit\'e de Paris}
\footnote{  \href{mailto:leonid.vainerman@unicaen.fr}{\textcolor{blue}{leonid.vainerman@unicaen.fr}}, \ \href{mailto:jean-michel.vallin@imj-prg.fr}{\textcolor{blue}{jean-michel.vallin@imj-prg.fr}}}

\begin{abstract} We characterize braided commutative Yetter-Drinfeld $C^*$-algebras over weak Hopf $C^*$-algebras
in categorical terms. Using this, we then study quotient type coideal subalgebras of a given weak Hopf $C^*$-algebra
$\mathfrak G$ and coideal subalgebras invariant with respect to the adjoint action of $\mathfrak G$. Finally, as an
example, we explicitly describe quotient type coideal subalgebras of the weak Hopf $C^*$-algebras associated with
Tambara-Yamagami categories.
\end{abstract}
\tableofcontents

\footnote {AMS Subject Classification [2010]{: Primary 18D10, Secondary
16T05, Tertiary 46L05}}

\footnote{ Keywords {: Coactions and corepresentations
of quantum groupoids, $C^*$-categories, reconstruction theorem.}}

\newpage
\newenvironment{dm}{\hspace*{0,15in} {{\bf Proof.}}}{$\square$}
%%%%%%%%%%%%%%%%%%%%%%%%%%%%%%%%%%%%%%%%%%%%%%%%%%%%%%%%%%%%%%%%%%%%%%%%%%%%%%%%%%%%%%%%%%%%%%%%%%%%%%%%%%%%%%
\begin{section}{Introduction}

This paper continues the study of coactions of weak Hopf C*-algebras on C*-algebras and their applications
which was initiated in \cite{VV1} and \cite{VV2}. Let us first recall our motivation.

It is known that any finite tensor category equipped with a fiber functor to the category of finite dimensional
vector spaces is equivalent to the representation category of some Hopf algebra - see, for example, \cite{EGNO},
Theorem 5.3.12. But many tensor categories do not admit a fiber functor, so they cannot be presented as representation
categories of Hopf algebras. On the other hand, T. Hayashi \cite{Ha} showed that any fusion category admits a tensor
functor to the category of bimodules over some semisimple (even commutative) algebra. Then it was proved in \cite{Ha},
\cite{Sz}, \cite{Os} that any fusion category is equivalent to the representation category of some algebraic structure
generalizing Hopf algebras called a {\it weak Hopf algebra} \cite{BNSz} or a {\it finite quantum groupoid} \cite{NV}.

The main difference between weak and usual Hopf algebra is that in the former the coproduct $\Delta$ is not
necessarily unital. In addition, a representation category of a weak Hopf algebra is, in general, multitensor,
i.e., its unit object is not necessarily simple (see, for example, \cite{EGNO}, 4.1). By this reason, in the present
paper we work mainly in the context of multitensor categories.

Apart from (multi)tensor categories, weak Hopf algebras have interesting applications to the subfactor theory. In particular,
for any finite index and finite depth $II_1$-subfactor $N\subset M$, there exists a weak Hopf $C^*$-algebra $\mathfrak G$
such that the corresponding Jones tower can be expressed in terms of crossed products of $N$ and $M$ with $\mathfrak G$
and its dual. Moreover, there is {\it a Galois correspondence} between intermediate subfactors in this Jones tower and
coideal $C^*$-subalgebras of $\mathfrak G$ - see \cite{NV2}. This motivates the study of coideal $C^*$-subalgebras
of weak Hopf $C^*$-algebras (in what follows - WHAs).

A unital $C^*$-algebra $A$ equipped with a {\it coaction} $\mathfrak a$ of a WHA $\mathfrak G=(B,\Delta, S,\varepsilon)$ is
called a $\mathfrak G$-$C^*$-algebra. When $A$ is a unital $C^*$-subalgebra of $B$ and $\mathfrak a=\Delta$, we call it a
coideal $C^*$-subalgebra or briefly a coideal of $B$.

%Note that if $\mathfrak G$
%is a usual Hopf $C^*$-algebra, then one can prove that necessarily $1_A=1_B$, so weak and usual coideals coincide.

The structure of the paper is as follows. Section 2 (Preliminaries) contains basic definitions and facts needed for the
comprehension of the main results of the paper. In particular, in Subsection 2.2 we describe three $C^*$-multitensor categories
associated with any weak Hopf $C^*$-algebra and in Subsection 2.3 we explain how to reconstruct a weak Hopf $C^*$-algebra
if one of these categories is given. Various results of this kind are known - see \cite{Sz}, \cite{Ha}, \cite{CE}, \cite{Pf},
\cite{Os}, and we present them in the form convenient for our goals.

It was shown in \cite{VV1} that any $\mathfrak G$-$C^*$-algebra $(A,\mathfrak a)$ corresponds to a pair $(\mathcal M,M)$,
where $\mathcal M$ is a module $C^*$-category with a generator $M$ over the category of unitary corepresentations of
$\mathfrak G$. Here we study an important special class of $\mathfrak G$-$C^*$-algebras - braided-commutative Yetter-Drinfel'd
$C^*$-algebras and characterize the corresponding $C^*$-module categories:
\begin{theorem} \label{catdual}
Given a WHA $\mathfrak G$, the following two categories are equivalent:

(i) The category $YD_{brc}(\mathfrak G)$ of unital braded-commutative Yetter-Drinfel'd $\mathfrak G$-$C^*$-algebras with unital $\mathfrak G$-
and $\hat{\mathfrak G}$-equivariant $*$-homomorphisms as morphisms.

(ii) The category $Tens(UCorep(\mathfrak G))$ of pairs $(\mathcal C,\mathcal E)$, where $\mathcal C$ is a $C^*$-multitensor category
whose associativities reduce to the changing of brackets and $\mathcal E:UCorep(\mathfrak G)\to\mathcal C$ is a unitary tensor functor
such that $\mathcal C$ is generated by the image of
$\mathcal E$. Morphisms $(\mathcal C,\mathcal E)\to(\mathcal C',\mathcal E')$ of this category are equivalence classes of pairs
$(\mathcal F,\eta)$, where $\mathcal F:\mathcal C\to\mathcal C'$ is a unitary tensor functor and $\eta:\mathcal F\mathcal E\to
\mathcal E'$ is a natural unitary monoidal functor isomorphism.

Moreover, given a morphism $[(\mathcal F,\eta)]:(\mathcal C,\mathcal E)\to(\mathcal C',\mathcal E')$, the corresponding homomorphism
of YD $\mathfrak G$-$C^*$-algebras is injective if and only if $\mathcal F$ is faithful, and it is surjective if and only if
$\mathcal F$ is full.
\end{theorem}
A similar result for compact quantum group coactions on C*-algebras was obtained earlier in \cite{NY1}. When it is possible, we
follow the same strategy. However, instead of tensor products over $\mathbb C$ we have to deal with tensor products over,
in general, non commutative algebras which makes many reasonings and calculations much more complicated.

%A general approach to the classification problem of indecomposable weak coideals of WHAs was developed in \cite{VV2}. Here,
In Section 4, we study, as an application of Theorem \ref{catdual}, coideals which belong to the category $YD_{brc}(\mathfrak G)$:
quotient type and invariant with respect to the adjoint action of a WHA and the relationship between them. We prove

%it is shown that any quotient type coideal is invariant. More difficult is the the inverse result:
\begin{theorem} \label{inv-qt}
Any quotient type coideal is invariant. Conversely, for any invariant coideal $I$ of $\mathfrak G$
%is covariantly isomorphic to a unique quotient type coideal, i.e.,
there exists a unique, up to isomorphism, quantum subgroupoid (i.e., a WHA $\mathfrak H$ equipped with an epimorphism
$\pi:\mathfrak G\to\mathfrak H$) such that $I$ is isomorphic as a $\mathfrak G$-$C^*$-algebra to the quotient type coideal
$I(\mathfrak H\backslash \mathfrak G)$.
\end{theorem}

Let us note that the coideal $(B,\Delta)$ is invariant (quotient type) if and only if $\mathfrak G$ is a usual Hopf algebra
and that invariant (quotient type) coideals form a sublattice of the lattice of all coideals of a WHA introduced in \cite{NV2}.

A concrete example illustrating the results of the paper is considered in Section 5. Namely, we describe invariant and
quotient type coideals of WHAs obtained by reconstruction from the Tambara-Yamagami categories \cite{TY} whose simple objects
are elements of a finite abelian group $G$ and one separate element $m$.
%$m$ satisfying the fusion rule $g\cdot h=gh,\ g\cdot m=m
%\cdot g=m,\ m^2=\underset{g\in G}\Sigma g,\ g^*=-g,\ m=m^*\ (g,h\in G)$. These categories are parameterized by non degenerate
%symmetric bicharacters $\chi:G\times G\to\mathbb C\backslash\{0\}$ and numbers $\tau=\pm |G|^{-1/2}$.
In particular, it is shown that a coideal is invariant if and only if it is of quotient type and that the lattice of invariant
(quotient type) coideals is isomorphic to the lattice of subgroups of $G$ completed by the new maximal element $G\sqcup m$.

%For any subset $K\subset G$, we shall denote $K^\perp:=\{g\in G|\chi(k,g)=1,\\forall k\in K\}.$
%The Hayashi's reconstruction theorem allows to construct a WHA  $\mathfrak G_{\mathcal T\mathcal Y}$ associated with
%$\mathcal T\mathcal Y(G,\chi,\tau)$. We recall this construction in Subsection 5.1.
%In fact, we give an explicit construction of representatives of all isomorphism classes of indecomposable finite dimensional
%$\mathfrak G_{\mathcal T\mathcal Y}$-$C^*$-algebras and indecomposable (weak) coideals of $\mathfrak G_{\mathcal T\mathcal Y}$.

{\em Notation:} for any category $\mathcal C$ we denote by $\Omega=Irr(\mathcal C)$ an exhaustive set of representatives of
the equivalence classes of its simple objects.

Our references are: to \cite{EGNO} for (multi)tensor categories, to \cite{Nes} for $C^*$-tensor categories and to \cite{NV}
for WHAs.

% Yetter-Drinfeld modules and the center construction for weak Hopf algebras were considered in \cite{Nen}.

\end{section}
%%%%%%%%%%%%%%%%%%%%%%%%%%%%%%%%%%%%%%%%%%%%%%%%%%%%%%%%%%%%%%%%%%%%%%%%%%%%%%%%%%%%%%%%%%%%%%%%%%%%%%%%%%%%%
\begin{section}{Preliminaries}

\begin{subsection}{Weak Hopf \texorpdfstring{$C^*$}{C*}-algebras} A {\it weak bialgebra} $\mathfrak G=(B,\Delta,\varepsilon)$
is a finite dimensional algebra $B$ with the comultiplication $\Delta : B\to B\otimes B$ and counit $\varepsilon
: B\to\mathbb C$ such that $(B,\Delta, \varepsilon)$ is a coalgebra and the
following axioms hold for all $b,c,d\in B$ :
\begin{enumerate}
\item[(1)] $\Delta$ is a (not necessarily unital) homomorphism : $\Delta(bc) = \Delta(b)\Delta(c).$

\item[(2)] The unit and counit satisfy the identities (we use the Sweedler leg notation $\Delta(c)=
c_{(1)}\otimes c_{(2)},\ (\Delta\otimes id_B)\Delta(c)=c_{(1)}\otimes c_{(2)}\otimes c_{(3)}$ etc.):
\begin{eqnarray*}
\varepsilon(bc_{(1)})\varepsilon(c_{(2)}d) &=&\varepsilon(bcd), \\
(\Delta(1)\otimes 1)(1\otimes \Delta(1)) &=& (\Delta\otimes id_B)\Delta(1).
\end{eqnarray*}
\end{enumerate}

A {\em weak Hopf algebra} is a weak bialgebra equipped with an antipode $S:B\to B$
which is an anti-algebra  and anti-coalgebra  homomorphism such that
\begin{eqnarray*}
m(id_B \otimes S)\Delta(b) &=& (\varepsilon\otimes id_B)(\Delta(1)(b\otimes 1)),\\
m(S\otimes id_B)\Delta(b) &=& (id_B \otimes \varepsilon)((1\otimes b)\Delta(1)),
\end{eqnarray*}
where $m$ denotes the multiplication.
%\end{enumerate}
%\medskip
%$\mathfrak
%G=(B,\Delta,S,\varepsilon)$ is a finite dimensional $C^*$-algebra $B$ with the
%comultiplication $\Delta : B\to B\otimes B$, counit $\varepsilon : B\to\mathbb C$, and
%antipode $S:B\to B$ such that $(B,\Delta, \varepsilon)$ is a coalgebra and the
%following axioms hold for all $b,c,d\in B$ :
%\begin{enumerate}
%\item[(1)] $\Delta$ is a (not necessarily unital)  $*$-homomorphism :
%$$
%\Delta(bc) = \Delta(b)\Delta(c), \quad \Delta(b^*) =
%\Delta(b)^*,
%$$
%\item[(2)] The unit and counit satisfy the identities (we use the Sweedler leg notation $\Delta(c)=
%c_1\otimes c_2,\ (\Delta\otimes id_B)\Delta(c)=c_1\otimes c_2\otimes c_3$ etc.):
%\begin{eqnarray*}
%\varepsilon(bc_1)\varepsilon(c_2d) &=&\varepsilon(bcd), \\
%(\Delta(1)\otimes 1)(1\otimes \Delta(1)) &=& (\Delta\otimes id_B)\Delta(1),
%\end{eqnarray*}
%\item[(3)]
%$S$ is an anti-algebra  and anti-coalgebra  map such that
%\begin{eqnarray*}
%m(id_B \otimes S)\Delta(b) &=& (\varepsilon\otimes id_B)(\Delta(1)(b\otimes 1)),\\
%m(S\otimes id_B)\Delta(b) &=& (id_B \otimes \varepsilon)((1\otimes b)\Delta(1)),
%\end{eqnarray*}
%where $m$ denotes the multiplication.
%\end{enumerate}

The right hand sides of two last formulas are called {\em target}
and {\em source counital maps} $\varepsilon_t$ and $\varepsilon_s$,
respectively. Their images are unital subalgebras of $B$ called
{\em target} and {\em source counital subalgebras} $B_t$ and $B_s$,
respectively. They commute elementwise, we have $S\circ\varepsilon_s =
\varepsilon_t\circ S$ and $S(B_t) =B_s$. We say that $B$ is {\em connected} if
$B_t \cap Z(B)= \mathbb{C}$ (where $Z(B)$ is the center of $B$), coconnected if
$B_t \cap B_s = \mathbb{C}$, and {\em biconnected} if both conditions are
satisfied.

Finally, if $B$ is a $C^*$-algebra and $\Delta(b^*) =\Delta(b)^*$, the
collection $\mathfrak G=(B,\Delta,S,\varepsilon)$ is called a {\em weak
Hopf $C^*$-algebra} (WHA). Then $B_t$ and $B_s$ are also $C^*$-subalgebras.

The dual vector space $\hat B$ has a natural structure of a WHA, namely 
$\hat{\mathfrak G}=(\hat B,\hat\Delta,\hat S,\hat\varepsilon)$
given by dualizing the structure operations of  $B$:
\begin{eqnarray*}
<\varphi\psi,\, b> &=& < \varphi\otimes\psi,\, \Delta(b)>, \\
<\hat\Delta(\varphi),\, b\otimes c> &=& <\varphi,\, bc>, \\
<\hat S(\varphi),\, b> &=& < \varphi,\, S(b)>, \\
< \phi^*,b> &=& \overline{ < \varphi,\, S(b)^*> },
\end{eqnarray*}
for all $b,c\in B$ and $\varphi,\psi\in \hat B$. The unit of $\hat B$ is
$\varepsilon$ and the counit is $1$.

The antipode $S$ is unique, invertible, and satisfies $(S\circ *)^2 =id_B$. Since
it was mentioned in \cite{Nik}, Remark 3.7 that problems regarding general WHAs can
be translated to problems regarding those with the property $S^2|_{B_t}=id$ which
are called {\em regular}, we will only consider such WHAs (see also \cite{Val6}).
In this case, there exists a canonical positive element $H$ in the center of $B_t$
such that $S^2$ is an inner automorphism implemented by $G=HS(H)^{-1}$, i.e.,
$S^2(b) = GbG^{-1}$ for all $b\in B$. The element $G$ is called the canonical
group-like element of $B$, it satisfies the relation\
$\Delta(G) =(G\otimes G)\Delta(1)= \Delta(1)(G\otimes G)$.

An element $\hat l\in\hat B$ is called a left integral (or a left invariant measure on $B$)
if $(id_B\otimes \hat l)\Delta = (\varepsilon_t\otimes \hat l)\Delta$. Similarly one gives
the definition of a right integral (or a right invariant measure on $B$). In any WHA there
is a unique positive left and right integral $h$ on $B$ such that $(id_B\otimes h)\Delta(1)=1$, 
called a {\em normalized Haar measure}.
%such that
%$$
%(id_B\otimes h)\Delta = (\varepsilon_t\otimes h)\Delta,\quad
%h\circ S =h,\quad h\circ \varepsilon_t = \varepsilon,\quad (id_B\otimes h)\Delta(1_B) = 1_B.
%$$
We will dehote by $H_h$ the GNS Hilbert space generated by $B$ and $h$ and by $\Lambda_h:
B\to H_h$ the corresponding GNS map.
\end{subsection}

%%%%%%%%%%%%%%%%%%%%%%%%%%%%%%%%%%%%%%%%%%%%%%%%%%%%%%%%%%%%%%%%%%%%%%%%%%%%%%%%%%%%%%%%%%%%%%%%%%%%%%%%%%%%%%%

\begin{subsection}{Three categories associated with a WHA}

{\bf 1. Unitary representations.}
Let $\mathfrak G=(B,\Delta, S,\varepsilon)$ be a weak bialgebra.
%We denote by $\varepsilon_t, \varepsilon_s$ its target and source counital maps, by $B_t$ and $B_s$ its target
%and source subalgebras, respectively, and by $G$ its canonical group-like element. We also denote by $h$ the
%normalized Haar measure of $\mathfrak G$.
Objects of the category $Rep(\mathfrak G)$ of representations of $\mathfrak G$ are finite rank left $B$-modules, simple
objects are irreducible $B$-modules and morphisms are $B$-linear maps. The tensor product of two objects $H_1,H_2
\in Rep(\mathfrak G)$ is the subspace $\Delta(1_B)\cdot(H_1\otimes H_2)$ of the usual tensor product together with the
action of $B$ given by $\Delta$. Tensor product of morphisms is the restriction of the usual tensor product of $B$-module
morphisms. Any $H\in Rep(\mathfrak G)$ is automatically a $B_t$-bimodule via $z\cdot v\cdot t:=zS(t)\cdot v,
\ \forall z,t\in B_t, v\in E$, and the above tensor product is in fact $\otimes_{B_t}$, moreover the $B_t$-bimodule
structure on $H_1\otimes_{B_t} H_2$ is given by $z\cdot \xi \cdot t =(z\otimes S(t))\cdot \xi,\ \forall z,t\in B_t, \xi
\in H_1\otimes_{B_t} H_2$. This tensor product is associative, so the associativity isomorphisms are trivial. The unit
object of $URep(\mathfrak G)$ is $B_t$ with the action of $B$ given by $b\cdot z:=\varepsilon _t(bz),\ \forall b\in B, z\in B_t$.

When $\mathfrak G$ is a WHA, it is natural to consider the category $URep(\mathfrak G)$ of its {\em unitary} representations
formed by finite rank left $B$-modules whose underlying vector spaces are Hilbert spaces $H$ with scalar product
$<\cdot,\cdot>$ satisfying $<b\cdot v,w>=<v,b^*\cdot w>,\quad\text{for all}\quad v,w\in H,\ b\in B.$ Then the above tensor
product is also a Hilbert space because $\Delta(1_B)$ is an orthogonal projection. The scalar product on $B_t$ is defined by
$<z,t>=h(t^*z)$.

For any morphism $f:H_1\to H_2$, let $f^*:H_2\to H_1$ be the adjoint linear map: $<f(v),w>=<v,f^*(w)>,\ \forall v\in H_1, w\in H_2$.
Clearly, $f^*$ is $B$-linear, $f^{**}=f$, $(f\otimes_{B_t} g)^*=f^*\otimes_{B_t} g^*$, and $End(H)$ is a $C^*$-algebra, for any
object $H$. So $URep(\mathfrak G)$ is a finite $C^*$-multitensor category (${\bf 1}$ can be decomposable).

The conjugate object for any $H\in URep(\mathfrak G)$ is the dual vector space $\hat H$ naturally identified
($v\mapsto \overline v$) with the conjugate Hilbert space $\overline H$ with the action of $B$ defined by $b\cdot\overline v=\overline{G^{1/2}S(b)^*G^{-1/2}\cdot v}$, where $G$ is the canonical group-like element of $\mathfrak G$. Then the
rigidity morphisms defined by
\begin{equation} \label{rigid}
R_H(1_B)=\Sigma_i (G^{1/2}\cdot \overline e_i\otimes_{B_t}\cdot e_i),\ \overline R_H(1_B)=\Sigma_i (e_i\otimes_{B_t}
G^{-1/2}\cdot\overline e_i),
\end{equation}
where $\{e_i\}_i$ is any orthogonal basis in $H$, satisfy all the needed properties - see \cite{BSz}, 3.6. Also,
it is known that the $B$-module $B_t$ is irreducible if and only if $B_t\cap Z(B)=\mathbb C 1_B$, i.e., if $\mathfrak G$ is
{\em connected}. So that, we have

\begin{proposition} \label{C*}
$URep(\mathfrak G)$ is a rigid finite $C^*$-multitensor category with trivial associativity
constraints. It is $C^*$-tensor if and only if $\mathfrak G$ is connected.
\end{proposition}

\begin{remark} \label{multi}
If $\{z_\alpha\}_{\alpha\in\Gamma}$ is the set of minimal orthoprojectors of $B_t\cap Z(B)$, then the trivial representation
denoted by ${\bf 1}$ admits a decomposition ${\bf 1}=\underset{\alpha\in\Gamma}\oplus {\bf 1}_\alpha$ with ${\bf 1}_\alpha$
irreducibles and according to \cite{EGNO}, Remark 4.3.4 we have:
\begin{equation} \label{decomp1}
URep(\mathfrak G)=\underset{\alpha,\beta\in\Gamma}\oplus \mathcal C_{\alpha\beta},
\end{equation}
where $\mathcal C_{\alpha\beta}$ are called the component subcategories of $\mathcal C$. Moreover:

(1) Every irreducible of $\mathcal C$ belongs to one of $\mathcal C_{\alpha\beta}$.

(ii) The tensor product maps $\mathcal C_{\alpha\beta}\times\mathcal C_{\gamma\delta}$ to $\mathcal C_{\alpha\delta}$
and equals to $0$ unless $\beta=\gamma$.

(ii) Every $\mathcal C_{\alpha\alpha}$ is a rigid finite $C^*$-tensor category with unit object ${\bf 1}_\alpha$.

(iv) The conjugate of any $X\in \mathcal C_{\alpha\beta}$ belongs to $\mathcal C_{\beta\alpha}$.
\end{remark}

%\begin{example} \label{G^K}
%Let $\mathfrak G_{TY}$ be the WHA constructed via the Hayashi's reconstruction functor from the Tambara-Yamagami category
%$\mathcal T\mathcal Y(G,\chi,\tau)$ and let $L$ be a subgroup of $G$. Then one can construct the sub WHA $\mathfrak G^L$ of
%$\mathfrak G_{TY}$ corresponding to $L$ using the reconstruction theorem  (see below) as follows. Let us equip the subcategory
%$\mathcal D=Vec_L$ with the tensor functor $\mathcal F:\mathcal
%D\to Corr_f(R)$, where $R=(\mathfrak G^L)_t$, which is the composition of the inclusion $i: \mathcal D\to \mathcal T\mathcal Y
%(G,\chi,\tau)$ and the Hayashi's canonical functor $\mathcal H$. Then one can check that the WHA structure of $\mathfrak G^L$ is
%obtained exactly like in \cite{M}, Theorem 1.3.4.

%The $C^*$-category $URep(\mathfrak G^L)$ is tensor and $UCorep(\mathfrak G^L)$ are multitensor but not tensor.
%Any object of $URep(\mathfrak G^L)$ and $UCorep(\mathfrak G^L)$ belongs automatically to $Corr_f(R)$, where $R=(\mathfrak G^L)_t$,
%this gives forgetful tensor functors.
%\end{example}

%\vskip 0.5cm
{\bf 2. Unitary comodules}

\begin{definition} \label{ucomod}
A  right unitary $\mathfrak G$-comodule is a pair  $(H,\mathfrak a)$, where $H$ is a Hilbert space with scalar product
$<\cdot,\cdot>$, $\mathfrak a : H \to H \otimes B$ is a bounded linear map between Hilbert spaces $H$ and $H\otimes H_h=
H\otimes \Lambda_h(B)$, and such that:

(i) $(\mathfrak a \otimes id_B)\mathfrak a= (id_H\otimes \Delta)\mathfrak a$;

(ii) $(id_H \otimes \varepsilon)\mathfrak a= id_H$;

(iii) $<v^{(1)},w>v^{(2)} =<v,w^{(1)}>S(w^{(2)})^*,\ \ \forall v,w\in H$,
where we used the leg notation $\mathfrak a(v)=v^{(1)}\otimes v^{(1)}$.

A morphism of unitary $\mathfrak G$-comodules $H_1$ and $H_2$ is a
linear map $T:H_1\to H_2$ such that $\mathfrak a_{H_2}\circ T=(T\otimes
id_B)\mathfrak a_{H_1}$ (i.e., a $B$-colinear map).

Right unitary $\mathfrak G$-comodules with {\bf finite dimensional} underlying Hilbert spaces and their morphisms
form a category which we denote by $UComod(\mathfrak G)$.

We say that two unitary $\mathfrak G$-comodules are equivalent (resp.,
unitarily equivalent) if the space of morphisms between them contains an
invertible (resp., unitary) operator.
\end{definition}

\begin{example} \label{ucoid}
Let us equip a right coideal $I\subset B$ with the scalar product
$<v,w>:=h(w^*v)$. Then the strong invariance of $h$ gives:
$$
<v^{(1)},w>v^{(2)}=(h\otimes id_B)((w^*\otimes 1_B)\Delta(v))=
$$
$$
=(h\otimes S^{-1})(\Delta(w^*)(v\otimes 1_B))=<v,w^{(1)}>S(w^{(2)})^*.
$$
\end{example}

%\begin{remark}
%\label{injective}
%By (ii) any coaction $\mathfrak a$ is injective.
%\end{remark}

If $(H,\mathfrak a)$ is a right  unitary $\mathfrak G$-comodule, then $H$ is naturally a unitary
left $\hat{\mathfrak G}$-module via
\begin{equation} \label{mult}
\hat b\cdot v:=v^{(1)}<\hat b,v^{(2)}>,\ \ \forall \hat b\in\hat B,\ v\in H.
\end{equation}
%The unitarity follows from the calculation
%$$
%<\hat b\cdot v,w> = <v^{(1)}<\hat b,v^{(2)}>,w> =  <\hat b,<v^{(1)},w>v^{(2)}> =
%$$
%$$
%= <\hat b,<v,w^{(1)}>S(w^{(2)})^*>=<v,w^{(1)}\overline{<\hat b,S(w^{(2)})^*>}>=<v,(\hat b)^*\cdot w>,
%$$ for all $v,w\in H$ and $\hat b\in\hat B$. In particular $H$ is a $\hat
%B_t$-bimodule.

Due to the canonical identifications $B_t\cong\hat B_s$ and $B_s\cong\hat
B_t$ given by the maps $z \mapsto\hat z=\varepsilon (\cdot z)$ and $t\mapsto
\hat t =\varepsilon (t\cdot)$, $H$ is also a $B_s$-bimodule via $z\cdot v\cdot
t=v^{(1)}\varepsilon (z v^{(2)} t)$, for all $z,t\in B_s,\ v\in V$. The maps
$\alpha,\beta:B_s\to B(H)$ defined by $\alpha(z)v:=z\cdot v$ and
$\beta(z)v:=v\cdot z$, for all $z\in B_s, v\in H$ are a $*$-algebra
homomorphism and antihomomorphism, respectively, with commuting images.
Indeed, for instance, for all $v,w\in H, z\in B_s$, one has:
$$
<\alpha(z)v,w>:=<v^{(1)}\varepsilon (zv^{(2)}),w>=\varepsilon (<v^{(1)},w>zv^{(2)})=
$$
$$
=\varepsilon (<v,w^{(1)}>zS(w^{(2)})^*)=<v,w^{(1)}>\overline{\varepsilon (S(w^{(2)})z^*)}=
$$
$$
=<v,w^{(1)}\varepsilon (S(z^*)w^{(2)})>=<v,\alpha(z^*)w^{(1)}\varepsilon (w^{(2)})>=<v,\alpha(z^*)w>.
$$
So that, $\alpha(z)^*=\alpha(z^*)$, and similarly for the map $\beta$.
%We have the following useful relations:
%\begin{equation} \label{compatible}
%\mathfrak a(\alpha(x)\beta(y)v)=v^{(1)}\otimes xv^{(2)}y\ \forall v\in H, x,y\in B_s.
%\end{equation}
%and
%\begin{equation} \label{compatible'}
%\alpha(x)\beta(y)v^{(1)}\otimes v^{(2)}=v^{(1)}\otimes S(x)v^{(2)}S(y)\ \forall v\in H, x,y\in B_s.
%\end{equation}

The  correspondence (\ref{mult}) is bijective since one has the inverse formula:
if $(b_i)_i$ is a basis for $B$ and $(\hat b_i)$ is its dual basis in $\hat B$,
then set:
\begin{equation}
\label{mult'}
\mathfrak a(v) = \underset{i} \sum (\hat b_i \cdot v) \otimes b_i \ \ \forall v \in H.
\end{equation}
Moreover, formulas (\ref{mult}) and (\ref{mult'}) also lead to a bijection of
morphisms, and we have two functors, $\mathcal F_1:UComod(\mathfrak G)\to
URep(\hat{\mathfrak G})$ and $\mathcal G_1:URep(\hat{\mathfrak G})\to
UComod(\mathfrak G)$, which are mutually inverse to each other. Hence, these categories are
isomorphic and we can transport various additional structures from $URep(\hat{\mathfrak G})$
to $UComod(\mathfrak G)$ and vice versa.

For instance, let us define tensor product of two unitary $\mathfrak G$-comodules,
$(H_1,\mathfrak a_{H_1})$ and $(H_2,\mathfrak a_{H_2})$. As a vector space, it is
$$
H_1\otimes_{\hat{B_t}} H_2:=\hat\Delta(\hat{1})(H_1\otimes H_2)={\hat{1}}_{(1)}\cdot H_1\otimes {\hat{1}}_{(2)}\cdot H_2
$$
and
%is generated by the elements $ x\otimes_{\hat B_t} y := \hat \Delta (\hat 1)\cdot(x \otimes y)$, where $x\in H_1,
%y\in H_2$, so it
can be identified with $ H_1\otimes_{{B_s}} H_2$ (see \cite{Pf1}, 2.2 or \cite{Ni}, Chapter 4).
The unitary comodule structure on $ H_1\otimes_{{B_s}} H_2$ is given by
$$
v \otimes_{B_s} w\mapsto v^{(1)}\otimes_{B_s} w^{(1)}\otimes v^{(2)}w^{(2)},\ \forall v\in H_1, w\in H_2.
$$

%\begin{lemma} \label{complexe}
%If $(H_1,\mathfrak a),\ (H_2,\mathfrak b)\in UComod(\mathfrak G)$, then the
%projection $P: H_1 \otimes H_2 \to H_1\otimes_{\hat{B_t}} H_2$ defined by $P(v) = \hat \Delta(\hat 1)\cdot v$, for all
%$v \in H_1 \otimes H_2$, satisfies
%$$
%P(x\otimes y) = x^{(1)}\otimes y^{(1)}\varepsilon (x^{(2)}y^{(2)}),\quad\text{for\ all}\ x \in H_1,\ y \in H_2.
%$$
%\end{lemma}
%The proof is the direct calculation using the axiom (2) of a weak Hopf algebra:
%$$
%{\hat{1}}_1\cdot x\otimes {\hat{1}}_2\cdot y=(x^1\otimes y^1)\varepsilon (x^2 1_1)\varepsilon (1_2 y^2)=
%$$
%$$
%=(x^1\otimes y^1)\varepsilon (x^2y^2).
%$$
%\begin{corollary} The linear map $\mathfrak a \otimes_{{B_s}} \mathfrak  b$ given by:
%$$
%v \otimes_{B_s} w\mapsto v^1\otimes_{B_s} w^1\otimes v^2w^2,\ \forall v\in H_1, w\in H_2,
%$$
%is a coaction of $\mathfrak G$ on $H_1\otimes_{{B_s}} H_2$ (i.e., satisfies Definition \ref{ucomod}, (i), (ii)).
%\end{corollary}

%The direct calculation shows that the tensor product coaction is unitary.

Thus, $UComod(\mathfrak G)$ is a multitensor category with  trivial associativity isomorphisms whose unit object
$(B_s,\Delta|_{B_s})$ is simple if and only if $\mathfrak G$ is coconnected.
%The left and right unit isomorphisms are:
%\begin{equation} \label{unit1}
%l_H : B_s \otimes_{B_s} H \to H, \ \ z \otimes_{B_s} v \mapsto z\cdot v, \ r_H : H \otimes_{B_s} B_s\to  H,
%\ \ v \otimes_{B_s} z \mapsto v\cdot z.
%\end{equation}
%One can check that these isomorphisms are unitary and their inverses are:
%\begin{equation} \label{unitinv}
%l_H^{-1}(v)=1_1\otimes_{B_s} v^1\varepsilon (1_2v^2)\quad\text{and}\quad r_H^{-1}(v)=v^1\otimes_{B_s}\varepsilon _s(v^2).
%\end{equation}
The conjugate object for $(H,\mathfrak a)\in UComod(\mathfrak G)$ is $(\overline H,\tilde{\mathfrak a})$ with
\begin{equation} \label{conj}
\tilde{\mathfrak a}(\overline v)=\overline{v^{(1)}}\otimes [{\hat G}^{-1/2}\rightharpoonup (v^{(2)})^*\leftharpoonup{\hat G}^{1/2}],
\end{equation}
where ${\hat b}\rightharpoonup b:=<\hat b,b_{(2)}b_{(1)}>,\ b\leftharpoonup{\hat b}:=<\hat b,b_{(1)}>b_{(2)}\ (\forall
b\in B,{\hat b}\in\hat B)$ are the Sweedler arrows and $\hat G$ is the canonical group-like element of $\hat{\mathfrak G}$.

%. The corresponding Hilbert space is
%$\overline H$. In what follows, we use the Sweedler arrows ${\hat b}\rightharpoonup b:=<\hat b,b_{(2)}b_{(1)}>,\
%b\leftharpoonup{\hat b}:=<\hat b,b_{(1)}>b_{(2)},\ \forall b\in B,{\hat b}\in\hat B$.

%\begin{lemma} \label{harpoon}
%The conjugate object for $(H,\mathfrak a)$ in $UComod(\mathfrak G)$ is $(\overline H,\tilde{\mathfrak a})$, where
%$$
%\tilde{\mathfrak a}(\overline v)=\overline{v^{(1)}}\otimes [{\hat G}^{-1/2}\rightharpoonup (v^{(2)})^*\leftharpoonup{\hat G}^{1/2}],
%$$
%and $\hat G$ is the canonical group-like element of the dual quantum groupoid $\hat{\mathfrak G}$.
%\end{lemma}

The rigidity morphisms are given by (\ref{rigid}) with $B_t$ replaced by $B_s$. For any morphism $f$, $f^*$ is the conjugate
linear map of the corresponding Hilbert spaces, the colinearity of $f$ implies that $f^*$ is colinear. So that, we have

\begin{proposition} \label{C**}
$UComod(\mathfrak G)$ is a strict rigid finite $C^*$-multitensor category isomorphic to $URep(\hat{\mathfrak G})$. It is
$C^*$-tensor if and only if $\mathfrak G$ is coconnected (i.e., $B_t\cap B_s=\mathbb C 1_B$).
\end{proposition}

{\bf 3. Unitary corepresentations.}

\begin{definition} \label{corepresentation}
A {\bf right unitary corepresentation} $U$ of $\mathfrak G$ on a Hilbert space $H_U$ is a partial isometry $U\in B(H_U)\otimes B$
such that:

(i) $(id\otimes\Delta)(U)=U_{12}U_{13}$.

(ii) $(id\otimes\varepsilon)(U)=id$.

A morphism between two right corepresentations $U$ and $V$ is a bounded linear map $T\in B(H_U,H_V)$ such that
$(T \otimes 1_B)U= V(T\otimes 1_B)$. We denote by $UCorep(\mathfrak G)$ the category whose objects are unitary
corepresentations on {\bf finite dimensional} Hilbert spaces and above mentioned morphisms.
\end{definition}

Any $H_U$ is a unitary right $B$-comodule via $v\mapsto U(v\otimes 1_B)$. Conversely, given $(H,\mathfrak a)\in
UComod(\mathfrak G)$, one can construct $V\in UCorep(\mathfrak G)$ as follows:
$$
V(x\otimes \Lambda_h y):=x^{(1)}\otimes \Lambda_h(x^{(2)}y)),\ \text{for\ all}\ x\in H,\ y\in B.
$$
%One also has:
%$$
% V^*(x\otimes \Lambda_h y):=x^{(1)}\otimes \Lambda_h(S(x^{(2)})y)),\ \text{for\ all}\ x\in H,\ y\in B.
%$$
Hence, the categories $UComod(\mathfrak G)$ and $UCorep(\mathfrak G)$ are isomorphic. The tensor product $U\otimes V$
equals $U_{13}V_{23}$ and acts on $H_U\otimes_{B_s} H_V$, the conjugate object $\overline U$ is the unitary
corepresentation acting on $\overline H_U$ via $\overline U(\overline x\otimes\Lambda_h(y))=\overline x^{[1)}\otimes
\Lambda_h((\overline x^{[2)})^*y)$, where $\tilde{\mathfrak a}(\overline x)$ is given by (\ref{conj}), the unit object
$U_\varepsilon\in B(B_s)\otimes B$ is defined by $z\otimes b\mapsto \Delta(1_B)(1_B\otimes zb),\ \forall z\in B_s, b
\in B$, and the rigidity morphisms are given by (\ref{rigid}) with $B_t$ replaced by $B_s$. For any morphism $T$,
$T^*$ is the conjugate linear map of the corresponding Hilbert spaces. Thus, we have
\begin{proposition} \label{C***}
$UCorep(\mathfrak G)$ is a strict rigid finite $C^*$-multitensor category isomorphic to $UComod(\mathfrak G)$ and to
$URep(\hat{\mathfrak G})$. It is $C^*$-tensor if and only if $\mathfrak G$ is coconnected.
\end{proposition}

\begin{remark} \label{mcprop}
%We denote by $\Omega$ an exhaustive set of representatives of the equivalence classes of irreducibles in $UCorep(\mathfrak G)$.
1. Using the leg notation $U=U^{(1)}\otimes U^{(2)}$, we define, for any $\eta,\zeta\in H_U$, the {\em matrix coefficient}
$U_{\eta,\zeta}:=<U^{(1)}\zeta,\eta>U^{(2)}\in B$ of $U$. If $\{\zeta_i\}$ is an orthonormal basis in $H_U$, denote $U_{i,j}
:=U_{\zeta_i,\zeta_j}$. Then the formula
$$
U=\oplus_{i,j} m_{i,j}\otimes U_{i,j},\ \text{where}\ m_{i,j}\ \text{are\ the\ matrix\ units\ of}\ B(H_U)\ \text{in  \ basis}\ \{\zeta_i\},
$$
defines a corepresentation of $\mathfrak G$ if and only if for all $i,j=1,...,dim(H_U)$:
\begin{equation} \label{prop}
\Delta(U_{i,j})= \Sigma^{dim(H_U)}_{k=1} U_{i,k}\otimes U_{k,j} ,\quad \varepsilon (U_{i,j})=\delta_{i,j},\quad
U_{i,j}=S(U_{j,i})^*.
\end{equation}

2. We also have $(U\otimes V)_{i,j,k,l}=U_{i,j}V_{k,l}$ for all $i,j=1,...,dim(H_U),k,l=1,...,dim(H_V),\ U,V\in UCorep(\mathfrak G)$.

3. For $U\in UCorep(\mathfrak G)$, denote $B_U:=Span\{U_{i,j}|i,j=1,...,dim(H_U)\}.$ Then (\ref{prop}) implies:
$\Delta(B_U)\subset \Delta(1_B)(B_U\otimes B_U)$, $B_U=S(B_U)^*$, $B_{\overline U}=(B_U)^*$.

4. a) $B_{\oplus^p_{k=1} U_k}=span\{B_{U_1},...,B_{U_p}\}$ for any finite direct sum of unitary corepresentations.
In particular, $B=\oplus_{x\in \Omega} B_{U^x}$.
%denoting $B_{U^x}$ by $B^x$ for irreducibles $U^x\ (x\in\Omega)$, we have
%$B=\oplus_{x\in \Omega} B_{U^x}$.

b) Decomposition $U\otimes V=\oplus_z d_z U^z$ with multiplicities $d_z$ implies $B_UB_V\subset \oplus_z B_{U^z}$, where
$z$ parameterizes the irreducibles of the above decomposition.
\end{remark}
\end{subsection}

\begin{subsection}{Reconstruction theorems.}
\label{generalreconstruction}
{\bf 1.} Let $\mathcal C$ be a rigid finite $C^*$-multitensor category with unit object ${\bf 1}$ and let $\mathcal J$
be a unitary tensor functor (see \cite{NT}, Definition 2.3.1) from $\mathcal C$ to the $C^*$-multitensor category $Corr_f(R)$
of finite dimensional Hilbert $R$-bimodules ($R$-correspon- dences), where $R=\mathcal J({\bf 1})$ is a finite
dimensional $C^*$-algebra. A discussion of the category $Corr_f(R)$ can be found in \cite{VV1}, pp. 86,87.

Put $H_U:=\mathcal J(U)$, for all $U\in\mathcal C$, in particular, $H^x:=\mathcal J(x)$, for all $x\in\Omega=Irr(\mathcal C)$.
Let $J_{U,V}: H_U\underset {R}\otimes H_V\to H_{U\otimes V}$ be the natural isomorphisms defining the tensor structure of 
$\mathcal J$ and choose an orthonormal basis $\{v^x_y|y\in\Omega_x:=\{1,...,dim (H^x)\}\}$ in each $H^x$.

Let $U^*$ be the conjugate of $U\in\mathcal C$, $R_U:{\bf 1}\to U^*\otimes U$ and $\overline R_U:{\bf 1}
\to U\otimes U^*$ be the corresponding rigidity morphisms.
%satisfying
%\begin{equation} \label{conjugate}
%(\overline R_U^*\otimes id_U)a^{-1}_{U,U^*,U}(id_U\otimes R_U)=id_U,\
%(R_U^*\otimes id_{U^*})a^{-1}_{U^*,U,U^*}(id_{U^*}\otimes\overline R_U).
%=id_{U^*}
%\end{equation}
Then the conjugate of $H_U$ is $H_{U^*}$ with the rigidity morphisms $J^{-1}_{U^*,U}\circ\mathcal J(R_U)$ and
$J^{-1}_{U,U^*}\circ\mathcal J(\overline R_U)$.
% and one can check that
%The unicity of a conjugate gives $H^{x^*}\cong Hom_R(H^x, R)$
%with rigidity morphisms $\tilde R_x:R\to H^{x^*}\underset{R}\otimes H^x$ and $\overline{\tilde R_x}:R\to H^x
%\underset{R}\otimes H^{x^*}$ given by
%$$
%\tilde R_x: 1_R\mapsto \underset{y\in \Omega_x}\Sigma v^{x^*}_y\underset{R}\otimes v^x_y,\quad
%\overline{\tilde R_x}: 1_R\mapsto \underset{y\in \Omega_x}\Sigma v^{x}_y\underset{R}\otimes v^{x^*}_y,
%$$
%where $\{v^{x}_y\}$ and $\{v^{x^*}_y\}$ are dual bases in $H^x$ and $H^{x^*}$, respectively.
The properties of the rigidity morphisms imply that the duality $<v,w>:=tr_R\circ\mathcal J(R^*_x)(J^{-1}_{x^*,x})^*
(v\underset{R}\otimes w)$, where $v\in H^{x^*}, w\in H^x$ and $tr_R$ is the trace of the left regular representation of $R$,
is non degenerate. Hence, there exist
isomorphisms $\Psi_x:\overline H^x\cong (H^x)^*\to H^{x^*}$ and $\Phi_x: H^x\to\overline H^{x^*}\cong (H^{x^*})^*$.

%$U^*$ is a left and right dual of $U$ in usual sense with $ev_U=R^*_U,\ coev_U=\overline R_U$, $ev'_U=\overline R^*_U,
%\ coev'_U=R_U$. We also have $(U^*)^*=U$ with $R_{U^*}=\overline R_U$, $\overline R_{U^*}=R_U$
%and $\overline{U\otimes V}=\overline V\otimes\overline U$.
%If $U$ is simple, the solution of (\ref{conjugate}) is unique up to $\lambda\in\mathbb C^{\times}$.

%If ${\bf 1}$ and $X$ are simple, there is a \alert{standard} solution of (1) satisfying
%$$
%R^*_X\circ R_X=\overline R^*_X\circ\overline  R_X=dim_q(X)\in End({\bf 1})=\mathbb C.
%$$
%$dim_q(X)>0$, it is uniquely determined by $X$ and called {\bf quantum dimension} of $X$. It coincides with right and left
%quantum dimensions in the usual sense, so $\mathcal C$ is spherical.

%A \alert{unitary} tensor functor is a tensor functor $(\mathcal F,\ \mathcal F_0,\ J_{X,Y})$ in the usual sense, where
%$\mathcal F:\mathcal C\to \mathcal C'$ is a $C^*$-functor, $\mathcal F_0:{\bf 1}_{\mathcal C'}\to {\bf 1}_{\mathcal C}$ and
%$J_{X,Y}:\mathcal F(X)\otimes \mathcal F(Y)\to \mathcal F(X\otimes Y)$ are unitary isomorphisms.

Next is a combined $C^*$-version of several reconstruction theorems scattered in various papers - see \cite{Sz}, \cite{Ha},
\cite{CE}, \cite{Pf}, \cite{Os}.

\begin{theorem} \label{genreconstruction}
A couple $(\mathcal C,\mathcal J)$ defines on the vector spaces
\begin{equation} \label{algebra}
B =  \underset {x \in \Omega} \bigoplus\ H^x\otimes\overline  H^x\quad\text{and}\quad
\hat B =  \underset {x\in\Omega}\bigoplus\ B(H^x)
\end{equation}
two  WHA structures, $\mathfrak G$ and $\hat{\mathfrak G}$, respectively,
dual to each other with respect to the bracket
$$
<A,w\otimes\overline v>=<Av,w>_x\quad\text{where}\quad x\in\Omega, A\in B(H^x),v,w\in H^x,
$$
%two  WHA structures, $\mathfrak G$ and $\hat{\mathfrak G}$, respectively,
%whose counital subalgebras are isomorphic (resp., antiisomorphic) to $R$ and
such that $\mathcal C\cong UComod(\mathfrak G)\cong UCorep(\mathfrak G)\cong URep(\hat{\mathfrak G})$.
\end{theorem}
{\bf A sketch of the proof.} Clearly, $\hat B$ is a $C^*$-algebra with usual matrix multiplication, conjugation and
unit. Since all of $H^x$ are Hilbert $R$-bimodules, there are homomorphisms $t:R\to\hat B$ and $s:R^{op}\to\hat B$
defined by $t(r)v=r\cdot v$ and $s(r)v=v\cdot r$, respectively (here $r\in R, v\in\underset {x\in\Omega}\oplus\ H^x$).

The coalgebra structure in $B$ dual to the algebra structure in $\hat B$, is:
\begin{equation} \label{coproduct}
\Delta(w\otimes\overline v) =  \underset {y\in\Omega_x}\bigoplus (w\otimes\overline v^{x}_y)_x
\otimes(v^{x}_y\otimes\overline v)_x,
\end{equation}
\begin{equation} \label{counit}
\varepsilon (w\otimes\overline v) =<w,v>_x,\quad\text{where}\quad v,w\in H^x.
\end{equation}

%It is known (see \cite{Sz}, \cite{Ha}, \cite{CE}, \cite{Os}) that $\hat B$
%the dual vector space
%\begin{equation} \label{dualalg}
%\hat B =  \undFirst, the coalgebra structure in $B$ dual to the algebra structure in $\hat B$, is:
%%\begin{equation} \label{coproduct}
%\Delta(w\otimes\overline v) =  \underset {y\in\Omega_x}\bigoplus (w\otimes\overline v^{x}_y)_x
%\otimes(v^{x}_y\otimes\overline v)_x,
%£\end{equation}
%\begin{equation} \label{counit}
%\varepsilon (w\otimes\overline v) =<w,v>_x,\quad\text{where}\quad v,w\in H^x.
%\end{equation}erset {x\in\Omega}\bigoplus\ B(H^x),
%\end{equation}
%where the duality is given, for all $x\in\Omega, A\in B(H^x),v,w\in H^x$ by:
%$$
%<A,w\otimes\overline v>=<Av,w>_x,
%$$

Then, as in \cite{CE}, 2.3.2, define the coproduct $\hat{\Delta}:\hat B\to\hat B \otimes\hat B$:
\begin{equation} \label{copr}
\hat\Delta(b) :=\eta\circ (J^{-1}\cdot b\cdot J)\quad (b\in \hat B, x,y\in \Omega),
\end{equation}
where $J:=\underset{x,y\in\Omega}\oplus J_{x,y},\ \eta:B(H^x)\underset{R}\otimes B(H^y)\to B(H^x)\otimes B(H^y)$
is the canonical map defined by
%the symmetric separability element $\underset{i\in I}e_i\otimes e^i\in R\otimes R:
$\eta(a_x\underset{R}\otimes c_y)=\underset{i\in I}\Sigma(s(e_i)a_x \otimes t(e^i) c_y)$ (here $a,c\in
\hat B,\{e_i\}$ and $\{e^i\}$ are dual bases of $R$ with respect to the duality $<a,b>=tr_R(L(a)L(b))$,
where $a,b\in R, L(a),L(b)$ are the corresponding left multiplication operators.
%and $tr_R$ is the trace of the left regular representation of $R$.

%In particular, $\hat\Delta(\hat 1)=\underset{i\in I}\Sigma (s(e_i)\otimes t(e^i))$, so $\varepsilon_t(b)=
%\underset{i\in I}\Sigma <b_{\bf 1},e_i> t(e^i),\ \varepsilon_s(b)=\underset{i\in I}\Sigma <b_{\bf 1},e^i> s(e^i)$.
The multiplication in $B$ dual to $\hat\Delta$, is as follows:
\begin{equation} \label{product}
(w\otimes\overline v)_x\cdot (g\otimes\overline h)_z =  ({J_{x,z}(w\underset{R}\otimes g)}\otimes \overline{J_{x,z}
(v\underset{R}\otimes h)})_{x\otimes z}\in H^{(x\otimes z)}\otimes\overline{H^{(x\otimes z)}},
\end{equation}
where $v,w\in H^x, g,h\in H^z$, for all $x,z\in\Omega$. For any $b=\underset{x\in\Omega}\oplus b_x$, the component
$b_{\bf 1}\in B(H^{\bf 1})\cong B(R)$, so $b_{\bf 1}(1_R)$ can be viewed as an element of $R$.
%, recall that H^1=R$, so we can associate with any $b\in \hat B$ an element $\tilde b_1:=b_1(1_R)\in R$.
Then the triple $(\hat B,\hat\Delta,\hat\varepsilon)$ with $\hat\varepsilon(b)=tr_R(L(b_{\bf 1}(1_R)))$
is a weak bialgebra.

%Let $\Omega_{\bf 1}\subset\Omega$ be the set of (classes of) irreducibles ${\bf 1_\omega}$ such that ${\bf 1}=\underset
%{\omega\in\Omega_{\bf 1}}\oplus\ {\bf 1_\omega}$ - see \cite{EGNO}, Corollary 4.3.2. With this notation, the triple $(B,
%\Delta,\varepsilon)$ with the above mentioned operations and unit ${\bf 1}_B=
%(\underset {\omega\in\Omega_{\bf 1}, x\in\Omega_\omega}\Sigma\ v^\omega_x)\otimes\overline
%{(\underset{\omega'\in\Omega_{\bf 1}, y\in\Omega_{\omega'}}\Sigma\ v^{\omega'}_y)}$ defines a
%dual weak bialgebra.

The antipode in $\hat B$ is defined by $\hat S(b)_x:=(\Psi_{x^*}\circ i_{x^*})(b_{x^*})^*(i^{-1}_{x^*}\circ\Phi_x)$,
where $b\in\hat B, x\in \Omega,\ i_x: H^x\to \overline H^x$ is a canonical antilinear isomorphism. Dually:
\begin{equation} \label{antipode}
S(w\otimes \overline v)=v^\natural\otimes{\overline w}^\flat,
(w\otimes \overline v)^*=w^\natural\otimes{\overline v}^\flat\ (\forall v,w\in H^x),
\end{equation}
where $w^\natural=\Psi_x(\overline w), {\overline v}^\flat=\Phi_x(v)$.
Any $H^x$ is a unitary right $B$-comodule via
$$
\mathfrak a_x(v)= \underset {y\in\Omega_x}\Sigma v\otimes v^x_y \otimes \overline{v^x_y} ,\quad
\text{where}\quad v\in H_x.
$$
which gives the equivalence $\mathcal C\cong UComod(\hat{\mathfrak G})\cong UCorep(\hat{\mathfrak G})\cong URep(\mathfrak G)$.
\hfill$\square$
\vskip 0.5cm
{\bf 2.} Let $(\mathcal C,\mathcal J)$ be as in Theorem \ref{genreconstruction}, let $\mathcal C'$ be a rigid finite
$C^*$-multitensor category and let $\mathcal P:\mathcal C'\to \mathcal C$ be
%with unit object ${\bf 1'}$ and $\Omega'==Irr(\mathcal C')$. Let there exist
a unitary tensor functor. Then Theorem \ref{genreconstruction} shows the existence of two WHAs, $\mathfrak G=(B,\Delta,
\varepsilon,S)$ and $\mathfrak G'=(B',\Delta',\varepsilon',S')$, generated by the couples $(\mathcal C,\mathcal J)$ and
$(\mathcal C',\mathcal J'=\mathcal J\circ\mathcal P)$, respectively, such that $UCorep(\mathfrak G)\cong UComod(\mathfrak G)
\cong\mathcal C$ and $UCorep(\mathfrak G')\cong UComod(\mathfrak G')\cong\mathcal C'$. The following theorem reveals the
structure of the functor $\mathcal P$.

\begin{theorem} \label{mapreconstr} In the above conditions, there is a WHA morphism $p:\mathfrak G'\to \mathfrak G$ such that
$\mathcal P$ viewed as a functor $UComod(\mathfrak G')\to UComod(\mathfrak G)$ is of the form $(V,\rho'_V)\mapsto (V,(id\otimes p)
\rho'_V)$, where $(V,\rho'_V)\in UComod(\mathfrak G')$. $\mathcal P$ is faithful (resp., full) if and only if $p$ is
is injective (resp., surjective).
%$\mathcal P$ is injective (resp., surjective) if and only if so is $p$.
\end{theorem}

Let us comment on the proof. For any $(V,\rho'_V)\in UComod(\mathfrak G')$ the condition $\mathcal J'=\mathcal J\circ\mathcal P)$
implies that $\mathcal P(V,\rho'_V)=(V,\rho_V)\in UComod(\mathfrak G)$, where $\rho_V:V\to V\otimes B$ is a right coaction. In
particular, $(B',\Delta')\in UComod(\mathfrak G')$, so $\mathcal P(B',\Delta')=(B',\rho_{B'})\in UComod(\mathfrak G)$, where
$\rho_{B'}:B'\to B'\otimes B$ is a right coaction. Then the composition $p:=(\varepsilon'\otimes id_B)\circ\rho_{B'}: B'\to B$
is a linear map. Theorem 3.5 of \cite{Wakui} proves that $p$ is a weak bialgebra morphism and that $\mathcal P(V,\rho'_V)=
(V,(id\otimes p)\rho'_V)$ for any $(V,\rho'_V)\in UComod(\mathfrak G')$. Corollary 3.6 of \cite{Wakui} shows that $\mathcal P$ is
an equivalence if and only if $p$ is an isomorphism.

In our context the comodules $(B',\Delta')$ and $\mathcal P(B',\Delta')=(B',(id_{B'}\otimes p)\Delta')$ are unitary
which gives for all $b,c\in B'$:
$$
<b_{(1)},c>(b_{(2)})=<b,c_{(1)}>S'(c_{(2)})^*,
$$ 
$$
<b_{(1)},c>p(b_{(2)})==<b,c_{(1)}>p(S(c_{(2)})^*).
$$
For $b=1_{B'}$ this implies $S(p(c))^*=p(S'(c)^*)$, for all $c\in B'$.

%{\bf To do: to show that $S(p(c))=p(S'(c))$ and $p(c)^*=p(c^*)$}

Then, the conjugate object for $(B',\Delta')$ in $UComod(\mathfrak G')$ is $(\overline B',\tilde{\Delta'})$, where
$$
\tilde{\Delta'}(\overline b)=\overline{b_{(1)}}\otimes [{\hat G}^{-1/2}\rightharpoonup (b_{(2)})^*\leftharpoonup{\hat G}^{1/2}],
$$
and $\hat G$ is the canonical group-like element of the dual WHA $\hat{\mathfrak G'}$. Let us note that $G$ and $\hat G$
belong to $B_t B_s$, so $p$ just sends them respectively to the canonical group-like elements of $\mathfrak G$ and its dual.
The conjugate object for $\mathcal P(B',\Delta)=(B',(id_{B'}\otimes p)\Delta')$ in $UComod(\mathfrak G)$ is described by a
similar formula. As $\mathcal P$ respects the rigidity of the categories in question, we have:
$$
\overline{b_{(1)}}\otimes p[{\hat G}^{-1/2}\rightharpoonup (b_{(2)})^*\leftharpoonup{\hat G}^{1/2}]=
\overline{b_{(1)}}\otimes [(\hat{p(G)}^{-1/2})\rightharpoonup (p(b_{(2)}))^*\leftharpoonup\hat{p(G)}^{1/2}]
$$
which gives by the invertibility of $\hat G$: $\overline{b_{(1)}}\otimes p(b_{(2)}^*)=\overline{b_{(1)}}\otimes
(p(b_{(2)}))^*$. Applying $\varepsilon'$ to the first leg, we get $p(b)^*=p(b^*)$, then also $S(p(b))=p(S'(b))
(\forall b\in B')$.

For the WHA $p(\mathfrak G')$ we have $UComod(p(\mathfrak G'))=\mathcal P(UComod(\mathfrak G'))$.
The functor $\mathcal P$ splits into the composition of the full functor $UComod(\mathfrak G')\to$\newline
$\mathcal P(UComod(\mathfrak G'))$ and the inclusion $\mathcal P(UComod(\mathfrak G')\to UComod(\mathfrak G)$.
Respectively, $p$ splits into the composition of the surjective WHA homomorphism $\mathfrak G'\to p(\mathfrak G')$
and the inclusion $p(\mathfrak G')\to \mathfrak G$. Then $\mathcal P$ is full if and only if $UComod
(p(\mathfrak G'))=UComod(\mathfrak G)$ or if and only if $p(\mathfrak G')=\mathfrak G$.

%Also, $B'=Vec(V_{ij})$, where $V_{ij}$ are matrix coefficients of $(V,\rho'_V)\in UComod(\mathfrak G')$ in some basis
%$\{v_i\}\in V$, and $p(B')=Vec(p(V'_{ij}))$. Now, if $\mathcal P$ is surjective, any element of $UComod(\mathfrak G)$
%is of the form $(V,(id_V\otimes p)\rho'_{V'})$ for some $(V,\rho'_V)\in UComod(\mathfrak G')$, so its matrix coefficients
%are also $p(V'_{ij}$ in the same basis. This implies that $p$ is surjective. Vice versa, if $p$ is surjective, we have
%$B=Vec(p(V'_{ij}))$, so $\mathcal P$ is surjective.

Also, $\mathcal P$ is faithful if and only if the first functor in the above decomposition is an
equivalence which happens if and only if $p:\mathfrak G'\to p(\mathfrak G')$ is an isomorphism of WHAs
or if only if the map $p:\mathfrak G'\to \mathfrak G$ is injective.
\hfill$\square$
% The proof implies that $p$ is injective (resp., surjective) if and only if so is $\mathcal P$.

%The corresponding pairs of total $C^*$-algebras in duality, $(B,\hat B)$ and $(B',\hat B')$, are given by the formulas
%(\ref{algebra}).

%Let $\Omega_{\bf 1}\subset\Omega$ be the set of (classes of) irreducibles ${\bf 1_\omega}$ such that ${\bf 1}=\underset
%{\omega\in\Omega_{\bf 1}}\oplus\ {\bf 1_\omega}$ - see \cite{EGNO}, Corollary 4.3.2. Then $R=
%= \underset{\omega\in\Omega_0}\oplus\ B(H^\omega)$, where $H^\omega=\mathcal J({\bf 1_\omega})$, and $J_{U,{\bf 1}}=
%J_{{\bf 1},U}=id_{H_U}$ for all $U\in\mathcal C$ - see \cite{EGNO}, Remark 2.4.6.

\end{subsection}

%%%%%%%%%%%%%%%%%%%%%%%%%%%%%%%%%%%%%%%%%%%%%%%%%%%%%%%%%%%%%%%%%%%%%%%%%%%%%%%%%%%%%%%%%%%%%%%%%%%%%%%%%%%%%%%%%%%%%%%%%%%%%%%

\begin{subsection}{Coactions.}
\begin{definition}\label{sieste}
A right coaction of a WHA $\mathfrak G$ on a unital $*$-algebra $A$, is a $*$-homomorphism
$\mathfrak a: A \to A \otimes B$  such that:

%1) $(A, \mathfrak a)$ is a right $B$-comodule.

1) $(\mathfrak a \otimes i)\mathfrak a= (id_A\otimes \Delta)\mathfrak a$.

2) $(id_A\otimes\varepsilon )\mathfrak a=id_A.$

3) $\mathfrak a(1_A) \in A \otimes B_t$.

One also says that $(A,\mathfrak a)$ is a $\mathfrak G$-$*$-algebra.
\end{definition}
If $A$ is a $C^*$-algebra, then $\mathfrak a$ is automatically continuous, even an isometry.

There are $*$-homomorphism $\alpha:B_s\to A$ and $*$-antihomomorphism $\beta:B_s\to A$ with commuting images defined by
$\alpha(x)\beta(y):=(id_A\otimes\varepsilon)[(1_A\otimes x)\mathfrak a(1_A)(1_A\otimes y)]$, for all $x,y\in B_s$.
We also have $\mathfrak a(1_A)=(\alpha\otimes id_B)\Delta(1_B)$,
\begin{equation} \label{compatible}
\mathfrak a(\alpha(x)a\beta(y))=(1_A\otimes x)\mathfrak a(a)(1_A\otimes y),
\end{equation}
and
\begin{equation} \label{compatible'}
(\alpha(x)\otimes 1_B)\mathfrak a(a)(\beta(y)\otimes 1_B)=(1_A\otimes S(x))\mathfrak a(a)(1_A\otimes S(y)).
\end{equation}
The set $A^{\mathfrak a} = \{a \in A | \mathfrak a(a) = \mathfrak a(1_A) (a \otimes 1_B) \}$ is a unital $*$-subalgebra
of $A$ (it is a unital $C^*$-subalgebra of $A$ when $A$ is a $C^*$-algebra) commuting pointwise with $\alpha(B_s)$.
A coaction $\mathfrak a$ is called {\it ergodic} if $A^{\mathfrak a} = \mathbb C 1_A$.

\begin{definition}\label{indecomp}
A $\mathfrak G-C^*$-algebra $(A,\mathfrak a)$ is said to be indecomposable if it cannot be presented as
a direct sum of two $\mathfrak G-C^*$-algebras.
\end{definition}
It is easy to see that $(A,\mathfrak a)$ is indecomposable if and only if $Z(A)\cap A^{\mathfrak a}=\mathbb C 1_A$.
Clearly, any ergodic $\mathfrak G-C^*$-algebra is indecomposable.

For any $(U,H_U)\in UCorep(\mathfrak G)$, we define the {\it spectral subspace} of $A$ corresponding to $(U,H_U)$ by
$$
A_U:=\{a\in A|\mathfrak a(a)\in \mathfrak a(1_A)(A\otimes B_U)\}.
$$
Let us recall the properties of the spectral subspaces:

(i) All $A_U$ are closed.

(ii) $A=\oplus_{x\in \Omega} A_{U^x}$.

(iii) $A_{U^x}A_{U^y}\subset \oplus_z A_{U^z},$ where $z$ runs over the set of all irreducible direct
summands of $U^x\otimes U^y$.

(iv) $\mathfrak a(A_U)\subset \mathfrak a(1_A)(A_U\otimes B_U)$ and $A_{\overline U}=(A_U)^*$.

(v) $A_\varepsilon $ is a unital $C^*$-algebra.
\end{subsection}

%%%%%%%%%%%%%%%%%%%%%%%%%%%%%%%%%%%%%%%%%%%%%%%%%%%%%%%%%%%%%%%%%%%%%%%%%%%%%%%%%%%%%%%%%%%%%%%%%%%%%%%

\begin{subsection}{Categorical duality.}

Let us note that the usage of $C^*$-multitensor categories allows to get without much effort the following
slight generalization
%t rid of the condition of
%coconnexity of a WHA in the formulation of
of the main result of \cite{VV1}:
\begin{theorem} \label{main}
Given a WHA $\mathfrak G$, the following categories are equivalent:

(i) The category of unital $\mathfrak G$-$C^*$-algebras with unital $\mathfrak G$-equivariant $*$-homomorphisms as morphisms.

(ii) The category of pairs $(\mathcal M, M)$, where $\mathcal M$ is a left module $C^*$-category with trivial module
associativities over   $UCorep(\mathfrak G)$ and $M$ is a generator in $\mathcal M$, with
equivalence classes of unitary module functors respecting the prescribed generators as morphisms.
\end{theorem}
In particular, given a unital $\mathfrak G$-$C^*$-algebra $A$, one constructs the $C^*$-category $\mathcal M=\mathcal D_A$ of
finitely generated right Hilbert $A$-modules which are equivariant, that is, equipped with a compatible right coaction \cite{BS1}.
Any its object is automatically a $(B_s,A)$-bimodule, and the bifunctor $U\boxtimes X:= H_U\otimes_{B_s} X\in \mathcal D_A$,
for all $U\in UCorep(\mathfrak G)$ and $X\in\mathcal D_A$, turns $\mathcal D_A$ into a left module $C^*$-category over
$UCorep(\mathfrak G)$ with generator $A$ and trivial associativities.

Vice versa, if a pair $(\mathcal M, M)$ is given, the construction of a $\mathfrak G$-$C^*$-algebra $(A,\mathfrak a)$ contains
the following steps. First, denote by $R$ the unital $C^*$-algebra $End_{\mathcal M}(M)$ and consider the functor $F:\mathcal C
\to Corr(R)$
defined on the objects by $F(U)=Hom_{\mathcal M} (M, U\boxtimes M)\ \forall U\in\mathcal C$. Here $X=F(U)$ is a right $R$-module
via the composition of morphisms, a left $R$-module via $rX=(id\otimes r)X$, the $R$-valued inner product is given by $<X,Y>=X^*Y$,
the action of $F$ on morphisms is defined by $F(T)X=(T\otimes id)X$. The weak tensor structure of $F$ (in the sense of \cite{Nes1})
is given by $J_{X,Y}(X\otimes Y)=(id\otimes Y)X$, for all $X\in F(U), Y\in F(V), U,V\in UCorep(\mathfrak G)$.

Then consider two vector spaces:
\begin{equation} \label{bspace}
A=\underset{x\in\Omega}\bigoplus A_{U^x}:=\underset{x \in\Omega}\bigoplus( F(U^x)\otimes\overline{H^x})
\end{equation}
and
\begin{equation} \label{recalgebra2}
\tilde{A}=\underset{U\in \|UCorep(\mathfrak G)\|}\bigoplus A_{U}:=\underset{U\in \|UCorep(\mathfrak G)\|}
\bigoplus(F(U)\otimes\overline{H_U}),
\end{equation}
where $F(U)=\underset{i}\bigoplus F(U_i)$ corresponds to the decomposition $U=\bigoplus U_i$ into irreducibles, and
$\|UCorep(\mathfrak G)\|$ is an exhaustive set of representatives of the equivalence classes of objects in
$UCorep(\mathfrak G)$ (these classes constitute a countable set). $\tilde A$ is a unital associative algebra with
the product
$$
(X\otimes\overline\xi)(Y\otimes\overline\eta)=(id\otimes Y)X\otimes(\overline\xi\otimes_{B_s}\overline\eta),
\ \forall (X\otimes\overline\xi)\in A_U, (Y\otimes \overline\eta)\in A_V,
$$
and the unit
$$
1_{\tilde A}=id_M\otimes\overline{1_B}.
$$
Note that $(id\otimes Y)X=J_{X,Y}(X\otimes Y)\in F(U\otimes V)$. Then, for any $U\in UCorep(G)$, choose isometries
$w_i:H_i\to H_U$ defining the decomposition of $U$ into irreducibles, and construct the projection $p_A:\tilde A\to A$ by
\begin{equation} \label{p}
p_A(X\otimes\overline\xi)=\underset{i}\Sigma(F(w_i^*)X\otimes\overline{w_i^*\xi}),\ \forall (X\otimes\overline\xi)\in A_U,
\end{equation}
which does not depend on the choice of $w_i$. Then $A$ is a unital $*$-algebra with the product $x\cdot y:=
p(xy)$, for all $x,y\in A$ and the involution  $x^*:=p(x^\bullet)$, where $(X\otimes\overline\xi)^\bullet:=
(id\otimes X^*) F(\overline R_U)\otimes\overline{\hat G^{1/2}\xi}$, for all $\xi\in H_U, X\in F(U), U\in UCorep(\mathfrak G)$.
Here $\overline R_U$ is the rigidity morphism from (\ref{rigid}). Finally, the map
\begin{equation} \label{coact}
\mathfrak a(X\otimes\overline\xi_i)=\underset{j}\Sigma(X\otimes\overline\xi_j)\otimes U^x_{j,i},
\end{equation}
where $\{\xi_i\}$ is an orthogonal basis in $H^x$ and $(U^x_{i,j})$ are the matrix elements of $U^x$ in this basis, is
a right coaction of $\mathfrak G$ on $A$. Moreover, $A$ admits a unique $C^*$-completion $\overline A$ such that $\mathfrak a$
extends to a continuous coaction of $\mathfrak G$ on it.

\begin{remark} \label{decomp}
We say that a $UCorep(\mathfrak G)$-module category is indecomposable if it is not equivalent to a direct sum of two nontrivial
$UCorep(\mathfrak G)$-module subcategories. Theorem \ref{main} implies that a $\mathfrak G-C^*$-algebra $(A,\mathfrak a)$ is
indecomposable if and only if the $UCorep(\mathfrak G)$-module category $\mathcal M$ is indecomposable.
\end{remark}
%\begin{remark} \label{ccoid}

%2) Let $I$ be a unital right coideal $*$-subalgebra of $B$. Then $I_{U^x}=I\cap B_{U^x}$ and $F(U^x)$ can be
%identified with a Hilbert subspace of $H^x\ (\forall x\in\Omega)$ and the coaction is the restriction of $\Delta$.
%\end{remark}

%\begin{example} \label{regular}
%The  $C^*$-algebra $B$ with the coproduct $\Delta$ viewed as $\mathfrak G$-$C^*$-algebra, corresponds to the
%$UCorep(\mathfrak G)$-module $C^*$-category $\mathcal Corr_f(B_s)$ with generator $M=B_s$: for any element $U\in UCorep(\mathfrak G)$
%and $N\in Corr_f(B_s)$, one defines $U\boxtimes N:=F(U)\otimes_{B_s} N$, where the functor $F: UCorep(\mathfrak G)\to Corr_f(B_s)\
%(F(U)=H_U)$ is the forgetful functor. Indeed, identifying $\mathcal M(B_s,H_U)$ with $H_U$, we get an isomorphism of the algebra
%$\tilde A$ constructed from the pair $(\mathcal M,B_s)$ onto $\tilde B=\underset{U}\bigoplus (H_U\otimes\overline H_U)$ and then an
%isomorphism $A\cong B=\underset{x\in\hat G}\bigoplus (H_x\otimes\overline H_x)$ such that $p_A:\tilde A\to A$ turns into the map
%$p_B:\tilde B\to B$ sending $\xi\otimes\overline\eta\in H_U\otimes\overline H_U$ nto the matrix coefficient $U_{\xi,\eta}$.
%\end{example}

\begin{remark} \label{image}
The equivalence between $\mathcal M$ and $\mathcal D_A$ maps any morphism $f:Hom_{\mathcal M}(U\otimes_{B_s} M, V\otimes_{B_s} M)$
to a morphism $\tilde f: H_U\otimes_{B_s} A\to H_V\otimes_{B_s} A\ (U,V\in UCorep(\mathfrak G)$). $\tilde f$ is an $A$-linear map
on the right intertwining $\delta_{H_U\otimes_{B_s} A}=U_{13}(id\otimes_{B_s}\alpha)$ and $\delta_{H_V\otimes_{B_s} A}=
V_{13}(id\otimes_{B_s}\alpha)$, so it can be written as
$$
\tilde f=\underset{i}\Sigma s_i\otimes_{B_s} a_i\in B(H_U,H_V)\otimes_{B_s} A
$$
acting by $\tilde f(\xi\otimes_{B_s} a)=\Sigma s_i(\xi)\otimes_{B_s} a_i a$, where $\xi\in H_U, a\in A$, and such that
$V_{13}(id\otimes\alpha)\tilde f=(\tilde f\otimes id)U_{13}(id\otimes\alpha)$.
%Moreover, $A$ is spanned by all such $a_i$'s.
\end{remark}
\end{subsection}
\end{section}
%%%%%%%%%%%%%%%%%%%%%%%%%%%%%%%%%%%%%%%%%%%%%%%%%%%%%%%%%%%%%%%%%%%%%%%%%%%%%%%%%%%%%%%%%%%%%%%%%%%%%%%%%%
\begin{section}{Yetter-Drinfel'd \texorpdfstring{$C^*$}{C*}-algebras over WHA}

\begin{subsection}{Basic definitions and results}

Let $\mathfrak G$ be a WHA, $\hat{\mathfrak G}$ be its dual and $(A,\mathfrak a)$ be a right unital $\mathfrak G$-$C^*$-algebra
which is also a left unital $\hat{\mathfrak G}$-$C^*$-algebra via a left coaction $\mathfrak b:A\to \hat B\otimes A$.
The coaction $\mathfrak b$ defines a right $B$-module algebra structure $\triangleleft: A\otimes B\to A$ by
$$
a\triangleleft b:=(b\otimes id_A)\mathfrak b(a),\quad\text{for\ all}\quad a\in A,\ b\in B.
$$
One can check that the following relations hold:
$$
a\triangleleft 1_B=a,\quad (ac)\triangleleft b=(a\triangleleft b_{(1)})(c\triangleleft b_{(2)})\ \forall a,c\in A,\ b\in B,
$$
\begin{equation} \label{triangle}
a^*\triangleleft b=(a\triangleleft S(b)^*)^*\quad \text{and}\quad 1_A\triangleleft b=1_A\triangleleft \varepsilon_s(b).
\end{equation}
Below we will use the leg notations for coactions and write $1$ instead of $1_B$.
%: for various comodule structures diffrent from $\Delta$ like $\mathfrak a,
%\delta$ etc. we denote $\mathfrak a(a)=a^{(1)}\otimes a^{(2)}$ etc., and $\Delta(b)=b_{(1)}\otimes b_{(2)}$,
%and instead of $1_B$ we write simply $1$.

\begin{lemma} \label{YDcoact}
The following two conditions are equivalent:

(i) the identity
\begin{equation} \label{yd}
\mathfrak a(a\triangleleft b)=(a^{(1)}\triangleleft b_{(2)})\otimes S(b_{(1)}) a^{(2)} b_{(3)},
\end{equation}
holds for all $a\in A,\ b\in B$.

(ii) the identity
\begin{equation} \label{yd1}
(id_{\hat B}\otimes\mathfrak a)\mathfrak b(a)=
W^*_{13}(\mathfrak b\otimes id_B)\mathfrak a(a)W_{13},
\end{equation}
holds for all $a\in A$, where the operator $W\in\mathcal L(\lambda_{h\otimes h} (B\otimes B))$
is defined by
$$
W(\lambda_{h\otimes h}(b\otimes c)):=\lambda_{h\otimes h}(\Delta(c)(b\otimes 1)),
$$
for all $b,c\in B$ (it is the adjoint of the {\it regular multiplicative
partial isometrty} $I$ of $B$ - see \cite{Val1}), and $W_{13}$ is the usual leg notation.
\end{lemma}

\begin{dm} As $A=\oplus_{x\in\Omega} A_{U^x}$, where $A_{U^x}$ is the spectral subspace of $A$ corresponding
to an irreducible corepresentation $U^x$ of $A$, it suffices to prove the statement for $a\in A_{U^x}$ only.
The matrix units $\{m^x_{i,j}\}$ of $B(H^x)$ with respect to some orthogonal basis $\{e_i\}$ in $H^x$ and
the corresponding matrix coefficients $U^x_{i,j}$ of $U^x$ with all possible $i,j,x$ form dual bases in $\hat B$
and $B$, respectively, so that $\mathfrak b$ can be restored from $\triangleright$ by
\begin{equation} \label{restore}
\mathfrak b(a)=\Sigma_{i,j}\ m^x_{i,j} \otimes (a\triangleleft U^x_{i,j}) ,\quad\text{for\ all}\ a\in A_{U^x}.
\end{equation}
Since $W=\oplus_{x\in\Omega}dim(H^x)U^x$ implements $\Delta$ and $\Delta(U^x_{i,j})=\Sigma_k U^x_{i,k}\otimes U^x_{k,j}$,
the right hand side of (\ref{yd1}) can be written for any $a\in A_{U^x}$ as
$$
(U^x_{13})^*(\mathfrak b\otimes id_B))\mathfrak a(a)U^x_{13}=
$$
$$
=\Sigma_{i,j,p,q} (m^x_{i,j}\otimes 1_A\otimes (U^x_{j,i})^*) (\mathfrak b(a^{(1)}\otimes a^{(2)}))(m^x_{q,p}\otimes 1_A
\otimes U^x_{q,p})=
$$
$$
=\Sigma_{i,j,p,q,r,s}(m^x_{i,j}m^x_{r,s}m^x_{q,p}\otimes(a^{(1)}\triangleleft U^x_{r,s})\otimes (U^x_{j,i})^*a^{(2)}U^x_{q,p})=
$$
$$
=\Sigma_{i,j,p,q} (m^x_{i,p} \otimes(a^{(1)}\triangleleft U^x_{j,q})\otimes (U^x_{j,i})^*a^{(2)}U^x_{q,p}).
$$
On the other hand, if (\ref{yd}) holds, the left hand side of (\ref{yd1}) can be written as
$$
(id_{\hat B}\otimes \mathfrak a)\mathfrak b(a)=\Sigma_{i,p} (m^x_{i,p}\otimes\mathfrak a(a\triangleleft U^x_{i,p})) =
$$
$$
=\Sigma_{i,p} (m^x_{i,p}\otimes(a^{(1)}\triangleleft (U^x_{i,p})_{(2)})\otimes S((U^x_{i,p})_{(1)}) a^{(2)} (U^x_{i,p})_{(3)})=
$$
$$
\Sigma_{i,p,j,q}(m^x_{i,p}\otimes(a^{(1)}\triangleleft U^x_{j,q})\otimes S(U^x_{i,j}) a^{(2)} U^x_{q,p}) =
$$
$$
\Sigma_{i,p,j,q}(m^x_{i,p}\otimes(a^{(1)}\triangleleft U^x_{j,q})\otimes (U^x_{j,i})^* a^{(2)} U^x_{q,p}).
$$
So (\ref{yd}) implies (\ref{yd1}). Conversely, writing in (\ref{yd1}) $\mathfrak b$ as above, we get (\ref{yd})
for any $b=U^x_{i,j}\ (x\in\Omega,i,j=1,...,dim H^x)$, $a\in A_{U^x}$ which gives the result.\hfill
\end{dm}
\begin{definition} \label{YD} (cf. \cite{NY1})
$A$ is a right-right Yetter-Drinfel'd (YD) $\mathfrak G$-$C^*$-algebra if
one of the above equivalent conditions is satisfied.

We say that a Yetter-Drinfel'd $\hat{\mathfrak G}$-$C^*$-algebra $A$ is
braided-commutative if
\begin{equation} \label{BC}
ab=b^{(1)}(a\triangleleft b^{(2)}),\quad\text{for\ all}\quad a,b\in A.
\end{equation}
\end{definition}

In particular, if $b\in A^\mathfrak a$, then $b^{(1)}\otimes b^{(2)}=1^{(1)}b\otimes 1^{(2)}$,
and since $b$ commutes with $1^{(1)}\in\alpha(B_s)$, the right hand side of (\ref{BC}) can be
written as $b1^{(1)}(a\triangleleft 1^{(2)})$. But (\ref{BC}) implies that
$1^{(1)}(a\triangleleft 1^{(2)})=a$, so $ab=ba$. Hence, $A^\mathfrak a\in Z(A)$.
%\vskip 0.5cm
%{\bf 2. Yetter-Drinfeld algebras and Drinfeld double}

Given a WHA $\mathfrak G$, let us construct a new WHA $D(\mathfrak G)$ called the {\bf Drinfel'd double} of $\mathfrak G$
as follows. The $C^*$-algebra of $D(\mathfrak G)$ is $B\otimes\hat B$, where $\hat{\mathfrak G}=(\hat B,\hat \Delta,\hat S,
\hat\varepsilon)$ is the dual of $\mathfrak G$. The coproduct $\Delta_D$ on $B\otimes\hat B$ is defined by
$$
\Delta_D=Ad(1\otimes\sigma\circ W\otimes 1_{\hat B})(\Delta\otimes\hat\Delta),
%(1\otimes\sigma\circ W^*\otimes 1_{\hat B}),
$$
where $W\in H_h\otimes H_{\hat h}$ is the multiplicative partial isometry
canonically associated with $\mathfrak G$ - see \cite{Val5}, and $\sigma$ is the flip. The antipode $S_D$
and the counit $\varepsilon_D$ on $B\otimes\hat B$  are defined, respectively, by
$$
S_D=Ad(W^*)(S\otimes\hat S)\quad\text{and}\quad\varepsilon_D=m(\varepsilon \otimes\hat\varepsilon).
$$
\begin{lemma} The collection $D(\mathfrak G)=(B\otimes\hat B,\Delta_D,S_D,\varepsilon_D)$ is a WHA.
\end{lemma}
\begin{dm} It suffices to note that the WHA constructed in \cite{NTV} and called there the
Drinfel'd double is dual to $D(\mathfrak G)$.\hfill\end{dm}

{\bf Theorem. $A$ YD $\mathfrak G$-$C^*$-algebra is the same as a $D(\mathfrak G)$-$C^*$-algebra.}

The proof is similar to the one of \cite{Nen}, Theorem 3.4.
\end{subsection}

%%%%%%%%%%%%%%%%%%%%%%%%%%%%%%%%%%%%%%%%%%%%%%%%%%%%%%%%%%%%%%%%%%%%%%%%%%%%%%%%%%%%%%%%%%%%%%%%%%%%%%%%%%%%%%%

\begin{subsection} {Categorical duality for Yetter-Drinfel'd algebras over WHAs}

Let us give the proof of Theorem \ref{catdual}. The condition that $\mathcal C$ is generated by $\mathcal E(\mathcal C)$
means that any object of $\mathcal C$ is isomorphic to a subobject of $\mathcal E(U)$ for some $U\in UCorep(\mathfrak G)$.
Assume without loss of generality that $\mathcal C$ is closed with respect to subobjects, but its unit object is not
necessarily simple.

Let us precise the equivalence relation on the set of pairs $(\mathcal F,\eta)$ in (ii). Given such a pair, we can consider,
for all $U,V\in UCorep(\mathfrak G)$, linear maps
$$
\mathcal C(\mathcal E(U),\mathcal E(V))\to \mathcal C'(\mathcal E'(U),\mathcal E'(V)):\quad T\mapsto \eta_V\mathcal F(T)\eta_U^{-1}.
$$
We say that two pairs, $(\mathcal F,\eta)$ and $(\mathcal F',\eta')$, are equivalent if the above maps are equal for all
$U,V\in UCorep(\mathfrak G)$.

The proof of Theorem \ref{catdual} will be done in several steps.

{\bf a) From YD $\mathfrak G$-$C^*$-algebras to $C^*$-multitensor categories.}

Given a braided commutative YD $\mathfrak G$-$C^*$-algebra $A$, let us show that the $C^*$-category $\mathcal D_A$
is in fact a $C^*$-multitensor category. We start with
\begin{remark} \label{delta} Recall the following relations:
%$\zeta\in H_U,\eta\in H_V\ U,V\in UCorep(\mathfrak G),a\in A$ :

1) $\delta_{H_U}(\zeta):=\zeta^{(1)}\otimes\zeta^{(2)}:=U(\zeta\otimes 1)$, where $\zeta\in H_U, U\in UCorep(\mathfrak G)$.

2) $\delta_{H_{U\otimes V}}(\zeta\otimes_{B_s}\eta):=U_{13}V_{23}(\zeta\otimes_{B_s}\eta\otimes 1)$ or $(\zeta\otimes_{B_s}\eta)^{(1)}=\zeta^{(1)}\otimes_{B_s}\eta^{(1)},\ (\zeta\otimes_{B_s}\eta)^{(2)}=\zeta^{(2)}\eta^{(2)}$,
where $\zeta\in H_U,\eta\in H_V,U,V\in UCorep(\mathfrak G)$.

3) $\delta_{H_U\otimes_{B_s} A}(\zeta\otimes_{B_s} a):=U_{13}(\zeta\otimes_{B_s}\mathfrak a(a))$ or $(\zeta\otimes_{B_s} a)^{(1)}=
\zeta^{(1)}\otimes_{B_s} a^{(1)},\ (\zeta\otimes_{B_s} a)^{(2)}=\zeta^{(2)} a^{(2)}$. Then $\delta_{H_{U\otimes V}\otimes_{B_s} A}
(\zeta\otimes_{B_s}\eta\otimes_{B_s} a)=\zeta^{(1)}\otimes_{B_s}\eta^{(1)}\otimes_{B_s} a^{(1)}\otimes\zeta^{(2)}\eta^{(2)} a^{(2)}$,
where $\zeta\in H_U,\eta\in H_V\ U,V\in UCorep(\mathfrak G),
a\in A$.

4) It follows from the equality $\mathfrak a(1_A)=(\alpha\otimes id)\Delta(1)$ (see \cite{VV1}) that
$(id_A\otimes\varepsilon_t)\mathfrak a(b)=\mathfrak a(1_A)(b\otimes 1)$. One can deduce from here, using 3)
and the relations $(id\otimes\varepsilon_t)U=(id\otimes\varepsilon_s)U=1$ that $(id\otimes\varepsilon_t)\delta_{H_U\otimes_A A}
(\zeta\otimes_A a)=\zeta\otimes_A a^{(1)}\otimes\varepsilon_t(a^{(2)})$.

\end{remark}

\begin{lemma} \label{2.2}

For any $X\in\mathcal D_A$, there exists a unique unital $*$-homomorphism $\pi_X:A\to \mathcal L_A(X)$
such that $\pi_X(a)(\zeta)=\zeta^{(1)}(a\triangleleft \zeta^{(2)})$ and $\delta_X(\pi_X(a)\zeta)=
(\pi_X\otimes id)\mathfrak a(a)\delta_X(\zeta)$, for all $a\in A$ and $\zeta\in X$.
\end{lemma}
\begin{dm} It suffices to consider $X=H^x\otimes_{B_s} A$ because $A$ is a generator of $\mathcal D_A$.
If $\{v^x_i\}$ is an orthonormal basis in $H^x$, Remark \ref{delta}, 2) is equivalent to

%call that the right comodule structure on $X=H^x\otimes_{B_s} A$ is given by
%$\delta_X(v\otimes_{B_s} b)=U^x_{13}(v\otimes_{B_s}\mathfrak a(b))$,
%where $v\in H^x, b\in A$, or
\begin{equation} \label{rightcom}
\delta_X(v^x_i\otimes_{B_s} b)=\underset{j}\Sigma[v^x_j\otimes_{B_s} b^{(1)}\otimes U^x_{j,i}b^{(2)}],
\end{equation}
here $U^x_{i,j}$ are the matrix coefficients of $U^x$.
%where $U^x_{i,j}$ are the matrix coefficients of $U^x$ and $m^x_{i,j}$ are the corresponding matrix units of $B(H_x)$.
The braided commutativity gives:
$$
(v^x_i\otimes_{B_s} b)^{(1)}(a\triangleleft (v^x_i\otimes_{B_s} b)^{(2)})=\underset{j}\Sigma(v^x_j\otimes_{B_s} b^{(1)})
(a\triangleleft (U^x_{j,i}b^{(2)}))=
$$
$$
=\underset{j}\Sigma(v^x_j\otimes_{B_s} b^{(1)}(a\triangleleft U^x_{j,i})\triangleleft b^{(2)})
=\underset{j}\Sigma(v^x_j\otimes_{B_s} (a\triangleleft U^x_{j,i}))b.
$$
Now it is clear how to define $\pi_X$ explicitly:
$$
\pi_X(a)=\Sigma_{i,j} m^x_{i,j}\otimes (a\triangleleft U^x_{j,i}),
$$
where $m^x_{i,j}$ are the corresponding matrix units of $B(H_x), a\in A$. This gives the first statement of the lemma.
In order to prove the second statement, take an arbitrary
$X\in\mathcal D_A$, then for any $a\in A$ and $\zeta\in X$ we have:
$$
\delta_X(\pi_X(a)\zeta)=\delta_X(\zeta^{(1)}(a\triangleleft \zeta^{(2)}))=(\zeta^{(1)}(a\triangleleft
\zeta^{(2)}))^{(1)}\otimes(\zeta^{(1)}(a\triangleleft \zeta^{(2)}))^{(2)}.
$$
The Yetter-Drinfel'd condition (\ref{yd}) shows that the last expression equals to
$$
\zeta^{(1)}(a\triangleleft \zeta^{(3)})^{(1)}\otimes \zeta^{(2)}(a\triangleleft \zeta^{(3)})^{(2)}=\zeta^{(1)}(a^{(1)}
\triangleleft \zeta^{(4)})
\otimes\zeta^{(2)} S(\zeta^{(3)})a^{(2)}\zeta^{(5)}=
$$
$$
=\zeta^{(1)}(a^{(1)}\triangleleft \zeta^{(3)})\otimes\varepsilon_t(\zeta^2)a^{(2)}\zeta^{(4)}.
$$
If again $X=H^x\otimes_{B_s} A$ and $\zeta=v^x_i\otimes_{B_s} b$, Remark \ref{delta}, 4) shows that
%$$
%(id\otimes\varepsilon_t)\delta_X(v^x_i\otimes_{B_s} b)=\underset{j}\Sigma(v^x_j\otimes_{B_s} b^1\otimes\varepsilon_t(U^x_{j,i}b^2))=
%$$
%$$
%=\underset{j}\Sigma(v^x_j\otimes_{B_s} b^1\otimes\varepsilon_t(U^x_{j,i}\varepsilon_t(b^2)))=
%\underset{j}\Sigma(v^x_j\otimes_{B_s} \alpha(1_1) b\otimes\varepsilon_t(U^x_{j,i}1_2))).
%$$
%because $\mathfrak a(b)=\mathfrak a(1_A)(b\otimes 1_B)$ since $\mathfrak a(1_A)=(\alpha\otimes id)\Delta(1_B)$ (see \cite{VV1}).
%Now, in any $C^*$-WHA we have the relation $1_1\otimes c 1_2 = \varepsilon_s((c_1)^*)\otimes (c_2)^*$ for all $c\in B$ and,
%for the matrix coefficients $U^x_{i,j}$, the relations $\varepsilon_t(U^x_{i,j})=\varepsilon_s(U^x_{i,j})=
%\delta_{i,j}1_B$, so we have
$(id\otimes \varepsilon_t)\delta_X(\zeta)=\zeta\otimes 1$ for the above $\zeta\in X$. This gives
$\delta_X(\pi_X(a)\zeta)=\zeta^{(1)}(a^{(1)}\triangleleft\zeta^{(2)})\otimes a^{(2)}\zeta^{(3)}$.
%$$
%(id\otimes\varepsilon_t)\delta_X(v^x_i\otimes_{B_s} b)=\underset{j}\Sigma(v^x_j\otimes_{B_s} b^1\otimes
%\varepsilon(u^x_{j,k}b^2)1_B=(v^x_i\otimes_{B_s} b)\otimes 1_B.
%$$
On the other hand,
$$
(\pi_X\otimes id)\mathfrak a(a)\delta_X(\zeta)=\pi_X(a^{(1)})\zeta^{(1)}\otimes a^{(2)}\zeta^{(2)}=\zeta^{(1)}(a^{(1)}
\triangleleft\zeta^{(2)})\otimes a^{(2)}\zeta^{(3)},
$$
%one checks all the needed properties of $\pi_X$ by direct calculations using formulas
%(\ref{restore}), (\ref{mod}), (\ref{BC}) and (\ref{yd}), similarly to the proof of \cite{NY1}, Lemma 2.2.
and we are done\hfill\end{dm}

This lemma implies that any $X\in\mathcal D_A$ is automatically a $\mathfrak G$-equivariant $(A,A)$-correspondence and
any $\mathfrak G$-equivariant endomorphism of the right Hilbert $A$-module $X$ is automatically an $(A,A)$-bimodule map.
Therefore, $\mathcal D_A$ is a full subcategory of the $C^*$-multitensor category of $\mathfrak G$-equivariant
$(A,A)$-correspondences. In order to show that $\mathcal D_A$ is invariant with respect to $\otimes_A$,
take $X,Y\in \mathcal D_A$ and prove two statements:

(i) $(X\otimes_A Y)\in \mathcal D_A$;

(ii) the left $A$-module structure on $X\otimes_A Y$ induced by that of $X$ is the same as the left $A$-module structure
given by Lemma \ref{2.2} using the coaction of $\mathfrak G$ and the right $A$-module structure on $X\otimes_A Y$.

The statement (ii) is proved by direct computations similar to those in the proof of Lemma \ref{2.2}.
In order to prove (i), it suffices to prove

\begin{lemma} \label{2.3} The map $T_{U,V}:X\otimes_A Y \to H_{U\otimes V} \otimes_{B_s} A$, where
$X=H_U\otimes_{B_s} A,\ Y=(H_V\otimes_{B_s} A)$ defined for all $\zeta\in H_U, \eta\in H_V, a,b\in A$ by
\begin{equation} \label{T}
T_{U,V}:(\zeta\otimes_{B_s} a)\otimes_A (\eta\otimes_{B_s} b)\mapsto \zeta\otimes_{B_s}\pi_Y(a)(\eta\otimes_{B_s} b),
\end{equation}
%where $\zeta\in H_V,\eta\in H_U,a,b\in A$ and $\zeta\to \zeta^1\otimes\zeta^2=U(\zeta\otimes 1_B)$ is a right $B$-comodule
%structure on $H_U$,
is a $\mathfrak G$-equivariant unitary isomorphism of right Hilbert $A$-modules and
\begin{equation} \label{tens}
T_{U\otimes V,W}(T_{U,V}\otimes_A id)=T_{U,V\otimes W}(id\otimes_A T_{V,W})\ (\forall\ U,V,W\in UCorep(\mathfrak G))
\end{equation}
\end{lemma}
\begin{dm} Clearly, (\ref{T}) defines $T_{U,V}$ as a right $A$-module isomorphism.
Let us note that the vectors of the form $(\zeta\otimes_{B_s} 1_A)\otimes_A (\eta\otimes_{B_s} 1_A)$
generate $X\otimes_A Y$ as a right $A$-module and that $T_{U,V}$ is isometric on these vectors. This implies that
$T_{U,V}$ is a unitary isomorphism of right Hilbert $A$-modules.

Let us check the $\mathfrak G$-equivariance of $T_{U,V}$, i.e., we must have the equality
$(T_{U,V}\otimes id_B)\delta_{X\otimes_A Y}=$ $\delta_{H_{U\otimes V} \otimes_{B_s} A}\circ T_{U,V}$.
Since $\pi_Y(a)(\eta\otimes_{B_s} b)=\eta^{(1)}\otimes_{B_s} (a\triangleleft \eta^{(2)})b$, we have:
$$
\delta_{X\otimes_A Y}((\zeta\otimes_{B_s} a)\otimes_A (\eta\otimes_{B_s} b))=(\zeta\otimes_{B_s} a)^{(1)}
\otimes_A(\eta\otimes_{B_s} b)^{(1)}\otimes
$$
$$
\otimes(\zeta\otimes_{B_s} a)^{(2)}(\eta\otimes_{B_s} b)^{(2)}=(\zeta^{(1)}\otimes_{B_s} b^{(1)})\otimes_A (\eta^{(1)}
\otimes_{B_s} a^{(1)})\otimes \zeta^{(2)} a^{(2)}\eta^{(2)}b^{(2)}
$$
Applying $T_{U,V}\otimes id_B$, we get $\zeta^{(1)}\otimes_{B_s} \eta^{(1)}\otimes_{B_s} (a^{(1)}\triangleleft
\eta^{(2)})b^{(1)}\otimes_{B_s}\zeta^{(2)}a^{(1)}\eta^{(3)}b^{(2)}$.

On the other hand,
$$
\delta_{H_{U\otimes V}\otimes_{B_s} A}\circ T_{U,V}[(\zeta\otimes_{B_s} a)\otimes_A (\eta\otimes_{B_s} b)]=
\delta_{H_{U\otimes V}\otimes_{B_s} A}(\zeta\otimes_{B_s} \eta^{(1)}\otimes_{B_s} (a\triangleleft \eta^{(2)})b)=
$$
$$
=\zeta^{(1)}\otimes_{B_s} \eta^{(1)}\otimes_{B_s} (a\triangleleft \eta^{(3)})^{(1)} b^{(1)}\otimes
\zeta^{(2)}\eta^{(2)}(a\triangleleft\eta^{(3)})^{(2)} b^{(2)}
$$
Applying (\ref{yd}), we see that the last expression equals to
$$
\zeta^{(1)}\otimes_{B_s} \eta^{(1)}\otimes_{B_s} (a\triangleleft \eta^{(4)}) b^{(1)}\otimes \zeta^{(2)}\eta^{(2)}
S(\eta^{(3)})a^{(2)}\eta^{(5)}b^{(2)}
$$
As $\eta^{(2)}S(\eta^{(3)})=\varepsilon_t(\eta^{(2)})$ and $\eta^{(1)}\otimes \varepsilon_t(\eta^{(2)})=1^{(1)}
\eta^{(1)})\otimes 1^{(2)}$, the last expression also equals to $\zeta^{(1)}\otimes_{B_s} \eta^{(1)}\otimes_{B_s}
(a^{(1)}\triangleleft \eta^{(2)})b^{(1)}\otimes_{B_s}\zeta^{(2)}a^{(1)}\eta^{(3)}b^{(2)}$.

Finally, the relation (\ref{tens}) can be justified by direct computations.
\hfill\end{dm}

\begin{corollary} \label{T_u,v} If $V=\underset{i,j}\Sigma m_{i,j}\otimes V_{i,j}$, then for all $\zeta\in H_U, \eta\in H_V, a,b\in A$
%The next formula follows from the proof of Lemma \ref{2.2}
$$
T_{U,V}(\zeta\otimes_{B_s} a)\otimes_A (\eta\otimes_{B_s} b)= \zeta\otimes_{B_s}\underset{i,j}\Sigma
m_{i,j}\eta\otimes_{B_s}(a\triangleleft V_{j,i}) b.
$$
\end{corollary}

Let us summarize the above mentioned results.
\begin{theorem} \label{2.4}
Let $A$ be a unital braided commutative YD $\mathfrak G-C^*$-algebra.
Then $\mathcal D_A$ is a $C^*$-multitensor category with tensor product $\otimes_A$ and trivial associativities equipped
with a unitary tensor functor $\mathcal E_A:UCorep(\mathfrak G)\to \mathcal D_A$ sending $U$ to $H_U\otimes_{B_s} A$ whose
structural unitary isomorphisms $T_{U,V}:\mathcal E_A(U)\otimes_{A} \mathcal E_A(V)\to \mathcal E_A(U\otimes V)$ are given
by (\ref{T}). Clearly, $\mathcal E_A(U_\varepsilon)=A={\bf 1}_{D_A}$.
\end{theorem}

{\bf b) From  $C^*$-multitensor categories to YD- $\mathfrak G$-$C^*$-algebras.}

Consider a pair $(\mathcal C,\mathcal E)\in Tens(UCorep(\mathfrak G))$. The category $\mathcal C$ is a left $UCorep(\mathfrak G)$-module
category: $U\boxtimes X:=\mathcal E(U)\otimes X,\ \forall X\in\mathcal C$, with generator $\mathcal E(\mathcal U_\varepsilon)=
{\bf 1}_{\mathcal C}$. So $R=End_{\mathcal C}({\bf 1}_{\mathcal C})$ and weak tensor functor $F$ sends any $U$ to $Hom_{\mathcal C}
({\bf 1}_{\mathcal C},\mathcal E(U))$. By Theorem \ref{main} we can construct a $\mathfrak G-C^*$-algebra $\overline A$ with right
coaction $\mathfrak a:\overline A\to \overline A\otimes B$. Now our goal is to prove

\begin{theorem} \label{tensYD}
The above $\mathfrak G-C^*$-algebra $\overline A$ has a natural structure of a unital braided commutative YD
$\mathfrak G-C^*$-algebra.
\end{theorem}

First, define a right $B$-module algebra structure on $A$ given by (\ref{bspace}). Let $\tilde A$ be an algebra
(\ref{recalgebra2}) with the projection $p_A:\tilde A\to A$ (\ref{p}) and $\tilde B=\underset{U\in\|UCorep(\mathfrak G)\|}
\bigoplus(H_U\otimes\overline H_U)$, $B=\underset{x\in\hat G} \bigoplus (H_x\otimes\overline  H_x)$ be the algebras with
the similar projection $p_B:\tilde B\to B$ (see \cite{VV1}, Example 6.7).
%mapping $\zeta\otimes\overline \xi\in H_U\otimes\overline H_U$ to $U_{\zeta,\xi}$,
Then define a linear map $\tilde{\triangleleft}:\tilde A\otimes B\to\tilde A$ for all $X\otimes\overline \eta\in A_U$,
$\zeta\otimes\overline \xi\in H_V\otimes\overline  H_V$ and $U,V\in UCorep(\mathfrak G)$:
\begin{equation} \label{module}
(X\otimes\overline\eta)\tilde{\triangleleft}(\zeta\otimes\overline\xi)=
(id\otimes X\otimes id)F(R_V)\otimes\overline{(\hat G^{-1/2}\cdot\overline\zeta\otimes_{B_s}\eta\otimes_{B_s}\xi)},
\end{equation}
where $R_V$ comes from (\ref{rigid}). Both sides of (\ref{module}) are in $A_{\overline V\otimes U\otimes V}$.
Identifying $B$ with the subspace of $\tilde B$, define a linear map $\triangleleft:A\otimes B\to A$ putting $a\triangleleft
b:=p_A(a\tilde{\triangleleft}b)$, for all $a\in A,b\in B$.

\begin{lemma} \label{2.6}
The map $\triangleleft$ defines a right $B$-module algebra structure on $A$ such that $p_A(a\tilde{\triangleleft } b)=
p_A(a)\triangleleft p_B(b)$, for all $a\in\tilde A, b\in\tilde B$.
\end{lemma}
\begin{dm}
Put $a=X\otimes\overline\eta\in F(U)\otimes\overline H_U$, $b=\zeta\otimes_{B_s}\overline\xi\in H_V\otimes_{B_s}\overline  H_V$
and choose isometries $u_i: H_{x_i}\to H_U$ and $v_j: H_{x_j}\to H_V$ defining the decompositions of $U$ and $V$ into irreducibles.
Then:
$$
p_A(a)\triangleleft p_B(b)=p_A(\underset{i,j}\Sigma(F(u^*_i) X\otimes\overline{u^*_i\eta})\tilde{\triangleleft }
(v^*_j \zeta\otimes_{B_s}\overline{v^*_j\xi}))=
$$
$$
=p_A(\underset{i,j}\Sigma(id\otimes F(u^*_i)X\otimes id)F(R_{V_{x_j}})\otimes\overline{(\hat G^{-1/2}\cdot\overline{v^*_j\zeta}\otimes_{B_s}u^*_i\eta\otimes_{B_s}
v^*_j\xi)}).
$$
On the other hand,
$$
p_A(a\tilde{\triangleleft } b)=p_A((id\otimes X\otimes id)F(R_V)\otimes\overline{(\hat G^{-1/2}\cdot\overline\zeta\otimes_{B_s}\eta\otimes_{B_s}\xi)})=
$$
$$
=p_A(\underset{i,j,k}\Sigma(F(\overline{v}^*_j)\otimes F(u^*_i)X\otimes F(v^*_k))F(R_V)\otimes\overline{(\overline{v^*_j
\hat G^{-1/2}\cdot\zeta}
\otimes_{B_s} u^*_i\eta\otimes_{B_s}v^*_k\xi)}),
$$
where the morphism $\overline{v}_j:H_{\overline{V}_{x_j}}=\overline H_{V_{x_j}}\to \overline H_V=H_{\overline V}$ is defined by $\overline{v}_j\overline\zeta=\overline{v_j\zeta}$. Since $v^*_j (\hat G^{-1/2}\cdot\xi)=\hat G^{-1/2}\cdot(v^*_j\xi),\ \forall j$,
$R_V=\underset{j}\Sigma(\overline v_j\otimes v_j)R_{V_j}$ and the partial isometries $v_j$ have mutually orthogonal images, the two
expressions are equal.

In order to show that $\triangleleft$ defines a right $B$-module algebra on $A$, take $a,b$ as above and $c=\mu\otimes_{B_s}
\overline\nu\in H_W\otimes_{B_s}\overline H_W$, where $W\in UCorep(\mathfrak G)$. Then:
$$
(a\tilde\triangleleft b)\tilde\triangleleft c=
(id\otimes X\otimes id)F(R_V)\otimes\overline{(\hat G^{-1/2}\cdot\overline\zeta\otimes_{B_s}\eta\otimes_{B_s}\xi)}\tilde\triangleleft
(\mu\otimes_{B_s}\overline\nu)=
$$
$$
(id\otimes id\otimes X\otimes id\otimes id)(id\otimes F(R_V)\otimes id)F(R_W)\otimes
$$
$$
\otimes\overline{(\hat G^{-1/2}\cdot\overline\mu
\otimes_{B_s}\hat G^{-1/2}\cdot\overline\zeta
\otimes_{B_s}\eta\otimes_{B_s}\xi\otimes_{B_s}\nu)}.
$$
The result belongs to $\tilde A_{\overline W\otimes(\overline V\otimes U\otimes V)\otimes W}$. On the other hand,
$$
a\tilde\triangleleft(bc)=(X\otimes\overline\eta)\tilde\triangleleft(\zeta\otimes_{B_s}\mu\otimes_{B_s}
\overline{\xi\otimes_{B_s}\nu})=
$$
$$
=(id\otimes X\otimes id)F(R_{V\otimes W})\otimes\overline{\hat G^{-1/2}
\cdot\overline{(\zeta\otimes_{B_s}\mu)}\otimes_{B_s}\eta\otimes_{B_s}(\xi\otimes_{B_s}\nu)}.
$$
This result belongs to $\tilde A_{\overline{V\otimes W}\otimes U\otimes V\otimes W}$ and is different from the previous one
because $\overline W\otimes\overline V\neq\overline{V\otimes W}$. But the map $\sigma:\overline H_W\otimes_{B_s}\overline H_V\to\overline{H_V\otimes_{B_s} H_W}$ defined by $\sigma(\overline\mu\otimes_{B_s}
\overline\zeta)=\overline{\zeta\otimes_{B_s}\mu}$ gives the equivalence of these corepresentations, so
$R_{V\otimes W}=(\sigma\otimes id\otimes id)(id\otimes R_V\otimes id)R_W$. Then, applying $p_A$ to the above elements,
we have an exact equality
$p_A((a\tilde\triangleleft b)\tilde\triangleleft c)=p_A(a\tilde\triangleleft (b c))$.

In order to check the relation $(ad)\triangleleft b=(a\triangleleft b_{(1)})(d\triangleleft b_{(2)})$, take $a=X\otimes\overline
\eta\in\tilde A_U, d=Y\otimes\overline\mu\in\tilde A_V, b=\zeta_i\otimes_{B_s}\overline\xi_j$, where $\{\zeta_i\otimes_{B_s}\overline
\xi_j\}_{i,j}$ is an orthonormal basis in $H_W\otimes_{B_s}\overline H_W$ and $U,V,W\in UCorep(\mathfrak G)$. Since $p_B(\zeta_i
\otimes_{B_s}\overline \xi_j)=W_{i,j}$ and
$\Delta(W_{i,j})=\underset{k}\Sigma(W_{i,k}\otimes W_{k,j})$, we have to show that
$$
p_A((ad)\tilde\triangleleft(\zeta_i\otimes_{B_s}\overline\xi_j))=\underset{k}\Sigma p_A((a\tilde\triangleleft(\zeta_i\otimes_{B_s}
\overline\xi_k))(d\tilde\triangleleft(\zeta_i\otimes_{B_s}\overline\xi_j))).
$$
The formula for the product in $\tilde A$ and (\ref{module}) give:
$$
p_A((ad)\tilde\triangleleft(\zeta_i\otimes_{B_s}\overline\xi_j))=
$$
$$=p_A((id\otimes X\otimes Y\otimes id)F(R_W)\otimes
\overline{\hat G^{-1/2}\cdot\overline\zeta_i\otimes_{B_s}\eta\otimes_{B_s}\mu\otimes_{B_s}\xi_j}).
$$
On the other hand,
$$
\underset{k}\Sigma (a\tilde\triangleleft(\zeta_i\otimes_{B_s}\overline\xi_k))
(d\tilde\triangleleft(\zeta_k\otimes_{B_s}\overline\xi_j)))=
$$
$$
=\underset{k}\Sigma((id\otimes X\otimes id)F(R_W)\otimes
\overline{\hat G^{-1/2}\cdot\overline{\zeta_i}\otimes_{B_s}\eta\otimes_{B_s}\xi_k})
$$
$$
((id\otimes Y\otimes id)F(R_W)\otimes\overline{\hat G^{-1/2}\cdot\overline{\zeta_k}\otimes_{B_s}\mu\otimes_{B_s}\xi_j})=
\underset{k}\Sigma(id\otimes X\otimes id\otimes id\otimes Y\otimes id)
$$
$$F(R_W\otimes R_W)\otimes\overline{(\hat G^{-1/2}\cdot\overline{\zeta_i}\otimes_{B_s}\eta\otimes_{B_s} \xi_k\otimes_{B_s}\hat G^{-1/2}\cdot\overline{\zeta_k}\otimes_{B_s}\mu\otimes_{B_s}\xi_j)}
$$
Since $\overline R_W({\bf 1})=\underset{k}\Sigma(\xi_k\otimes_{B_s}\hat G^{-1/2}\cdot\overline \zeta_k)$, and $\overline R_W$ is,
up to a scalar factor, an isometric embedding of ${\bf 1}$ to $W\otimes\overline W$, by applying $p_A$ to this element,
we get
$$
(id\otimes X\otimes F(\overline R^*_W)\otimes Y\otimes id)F(R_W\otimes R_W)\otimes\overline{(\hat G^{-1/2}\cdot\overline{\zeta_i}\otimes_{B_s}\eta\otimes_{B_s}\otimes_{B_s}\mu\otimes_{B_s}\xi_j)}
$$
Since $(\overline R^*_W\otimes id)(id\otimes R_W)=id_W$, this is equal to $p_A((ad)\tilde\triangleleft(\zeta_i\otimes_{B_s}
\overline\xi_j))$.
\hfill\end{dm}

Let us check now the compatibility of $\triangleleft$ with the involution.
\begin{lemma} \label{2.7} We have $a^*\triangleleft b=(a\triangleleft S(b)^*)^*$ for all $a\in A, b\in B$.
\end{lemma}
\begin{dm} Recall that if $a=X\otimes\overline\eta\in\tilde A_U$, then $a^\bullet=(id\otimes X^*)F(R_U)\otimes\overline{\overline{\hat G^{1/2}\cdot\eta}}\in
\tilde A_{\overline U}$. If $b=\zeta\otimes_{B_s}\overline\xi\in H_V\otimes_{B_s}\overline H_V$, put $b^\bullet=\xi\otimes_{B_s}
\overline\zeta$. Since $p_B(b^\bullet)=p_B(\xi\otimes_{B_s}\overline\zeta)=V_{\xi,\zeta}=(S(V_{\zeta,\xi}))^*=(S(p_B(b))^*$, we have
to prove that $p_A(a^*\tilde\triangleleft b)=p_A((a\tilde\triangleleft b^\bullet)^*)$.
Let us compute:
$$
(X\otimes\overline\eta)^\bullet\tilde\triangleleft(\zeta\otimes_{B_s}\overline\xi)=
$$
$$
=(id\otimes id\otimes X^*\otimes id)(id\otimes F(R_U)\otimes id)F(R_V)\otimes\overline{(\hat G^{-1/2}\cdot\overline\zeta \otimes_{B_s}
\overline{\hat G^{1/2}\cdot\eta}\otimes_{B_s}\xi)}.
$$
On the other hand,
$$
(a\tilde\triangleleft(\zeta\otimes_{B_s}\overline\xi)^\bullet)^\bullet=((id\otimes X\otimes id)F(R_V)\otimes\overline{(\hat G^{-1/2}\cdot\overline\xi\otimes_{B_s}\eta\otimes_{B_s}\zeta)})^\bullet=
$$
$$
=(id\otimes(id\otimes X\otimes id)F(R_V))^*F(R_{\overline V\otimes U\otimes V})\otimes
\overline{(\xi\otimes_{B_s}\overline{\hat G^{1/2}\cdot\eta}\otimes_{B_s}\overline{\hat G^{1/2}\cdot\zeta})}.
$$
Comparing these expressions and using the fact that $\hat G^{-1/2}\cdot\overline\xi=\overline{\hat G^{1/2}\cdot\xi}$, we see
that they are not equal only by the reason that the corepresentations $\overline V\otimes U\otimes V$ and $\overline{V\otimes
\overline U\otimes\overline V}$ are not equal. But they are equivalent via the map $\sigma(\overline\zeta\otimes_{B_s}\eta\otimes_{B_s}\xi)=\overline{\overline\xi\otimes_{B_s}
\overline\eta\otimes_{B_s}\zeta}$ which gives the relation
$$
R_{\overline V\otimes U\otimes V}
=(\sigma\otimes id\otimes id\otimes id)(id\otimes\overline R_V\otimes id\otimes id)
(id\otimes R_U\otimes id)R_V)
$$
Since $(R^*_V\otimes id)(id\otimes\overline R_V)=id_{\overline V}$, we have
$$
(id\otimes(id\otimes X\otimes id)F(R_V))^*F(R_{\overline V\otimes U\otimes V})=
$$
$$
=\sigma(id\otimes id\otimes X^*\otimes id)
(id\otimes F(R_U)\otimes id)F(R_V).
$$
Hence, the images of these expressions after applying $p_A$ are equal
\hfill\end{dm}

Now let us check the Yetter-Drinfel'd relation (\ref{yd}).

\begin{lemma} \label{2.8} For all $a\in A$ and $b\in B$ we have
$$
\mathfrak a(a\triangleleft b)=(a^{(1)}\triangleleft b_{[2)})\otimes S(b_{[1)}) a^{(2)} b_{(3)}
$$
\end{lemma}
\begin{dm} Let $U,V\in UCorep`(\mathfrak G)$ and $\{\eta_i\}\in H_U$, $\{\zeta_j\}\in H_V$ be two orthonormal
bases. For the simplicity, consider $\zeta_j$ as eigenvectors of the strictly positive operator $\zeta\mapsto
\hat G\cdot\zeta$ in $H_V:\hat G\cdot\zeta_j=\lambda_j(V)\zeta_j$. Then one has the following relations between the
matrix coefficients of $V$ with respect to $\{\zeta_j\}$ and of $\overline V$ with respect to $\{\overline\zeta_j\}$: $\lambda^{-1/2}_j(V)\lambda^{1/2}_k(V)\overline V_{j,k}= V^*_{j,k}=S(V_{k,j})$.

Now take $a=X\otimes\overline\eta_{k_0}\in\tilde A_U$ and $b=V_{i_0,j_0}$, then we have, using (\ref{coact}):
%$$
%\mathfrak a(a)=\underset{k}\Sigma(X\otimes\overline\xi_k)\otimes U_{k,k_0},
%$$
%from where
$$
(a^{(1)}\triangleleft b_{(2)})\otimes S(b_{(1)})a^{(2)}b_{(3)}=\underset{i,j,k}\Sigma[p_A(X\otimes\overline\eta_k)
\triangleleft V_{i,j}]\otimes S(V_{i_0,i})U_{k_0,k}V_{j_0,j}=
$$
$$
=\underset{i,j,k}\Sigma[p_A((id\otimes X\otimes id)F(R_V)\otimes\overline{(\overline\zeta_i\otimes_{B_s}\eta_k\otimes_{B_s}\zeta_j)}
\otimes \lambda^{1/2}_i(V) S(V_{i_0,i})U_{k,k_0}V_{j,j_0}.
$$
On the other hand,
$$
\mathfrak a(a\triangleleft b)=\mathfrak a(p_A((id\otimes X\otimes id)F(R_V) \otimes \overline{(\overline\zeta_{i_0}\otimes_{B_s}\eta_{k_0}\otimes_{B_s}\zeta_{j_0})}))=
$$
$$
=\lambda^{1/2}_{i_0}(V)\underset{i,j,k}\Sigma [p_A((id\otimes X\otimes id)F(R_V) \otimes \overline{(\overline\zeta_{i}\otimes_{B_s}\eta_{k}\otimes_{B_s}\zeta_{j})})\otimes\overline V_{i,i_0}U_{k,k_0}V_{j,j_0}.
$$
As $\lambda^{-1/2}_{i}(V)\lambda^{1/2}_{i_0}(V)\overline V_{i,i_0}=S(V_{i_0,i})$, the two expressions
are equal.
\hfill\end{dm}

Finally, let us check the braided commutativity relation (\ref{BC})

\begin{lemma} \label{2.9} For all $a,b\in A$ we have $ab=b^{(1)}(a\triangleleft b^{(2)})$.
\end{lemma}
\begin{dm} Let $a=p_A(X\otimes\overline\eta), b=p_A(Y\otimes\overline\zeta_i)$, where $X\in F(U), Y\in F(V), \eta\in H_U,
U,V\in UCorep(\mathfrak G)$ and bases $\{\eta_i\}\in H_U$ and $\{\zeta_j\}\in H_V$ as above. Then we compute:
$$
b^{(1)}(a\triangleleft b^{(2)})=\underset{j}\Sigma p_A((Y\otimes\overline\zeta_j)((X\otimes\overline\eta)\tilde
\triangleleft V_{j,i}))=
$$
$$
=\underset{j}\Sigma p_A((Y\otimes\overline\zeta_j)((id\otimes X\otimes id)F(R_V)\otimes\lambda^{1/2}_j\overline
{(\overline\zeta_j\otimes_{B_s}\eta\otimes_{B_s}\zeta_i)})=
$$
$$
=\underset{j}\Sigma p_A((id\otimes(id\otimes X\otimes id)F(R_V))Y\otimes\lambda^{1/2}_j(\overline\zeta_j\otimes\zeta_j
\otimes_{B_s}\overline\eta\otimes_{B_s}\overline\zeta_i))=
$$
$$
=p_A((id\otimes id\otimes X\otimes id)(id\otimes F(R_V))Y\otimes\overline{(\overline R_V(1_B)\otimes_{B_s}\eta\otimes_{B_s}\zeta_i)}
$$
Since $\overline R_V$ is, up to a scalar factor, an isometric embedding of ${\bf 1}$ into $V\otimes\overline V$, the last
expression equals
$$
p_A((F(\overline R^*_V)\otimes X\otimes id)(id\otimes F(R_V))Y\otimes\overline{(\eta\otimes_{B_s}\zeta_i)}=
p_A((id\otimes Y)X\otimes\overline{(\eta\otimes_{B_s}\zeta_i)}),
$$
which is exactly $ab$.\hfill\end{dm}
%Moreover, the map $\triangleleft$ satisfies all the relations (\ref{triangle})

%Then, the definition of the involution in $A$, the construction of the coaction (\ref{coact}), and the
%considerations similar to the ones in the proofs of \cite{NY1}, Lemmas 2.7, 2.8 and 2.9 allow to deduce all the relations
%(\ref{triangle}), (\ref{yd}) and (\ref{BC}).
Passing to the $C^*$-completion of $A$, we finish the proof of Theorem \ref{tensYD}.

%%%%%%%%%%%%%%%%%%%%%%%%%%%%%%%%%%%%%%%%%%%%%%%%%%%%%%%%%%%%%%%%%%%%%%%%%%%%%%%%%%%%%%%%%%%%%%%%%%%%%%%%%%%%%%%

{\bf c) Functoriality.} Given a morphism $\overline A_0\to\overline  A_1$ in $YD_{brc}(\mathfrak G)$, the map $X\to X
\otimes_{\overline A_0}\overline  A_1$ defines a unitary functor $\mathcal D_{\overline A_0}\to \mathcal D_{\overline A_1}$
(see \cite{VV1}, Theorem 4.12). By Theorem \ref{2.4}, both $\mathcal D_{\overline A_0}$ and $\mathcal D_{\overline A_1}$ are            $C^*$-multitensor categories, and similarly
to the proof of Lemma \ref{2.3} one shows that the isomorphisms
$$
(X\otimes_{\overline A_0}\overline A_1)\otimes_{\overline A_1}(Y\otimes_{\overline A_0} \overline A_1)\cong
(X\otimes_{\overline A_0} Y)\otimes_{\overline A_0} \overline A_1
$$
defined by $(x\otimes_{\overline A_0} a)\otimes_{\overline A_1} (y\otimes_{\overline A_0} b)\mapsto x\otimes_{\overline A_0} y^{(1)}
\otimes_{\overline A_1}(a\triangleleft (y^{(2)})b$, for all $x\in X\in \mathcal D_{\overline A_0}$
and $y\in Y\in \mathcal D_{\overline A_1}$, define a tensor structure on this functor. This functor together with obvious isomorphisms
$\eta_U:(H_U\otimes_{B_s} \overline A_0)\otimes_{\overline A_0} \overline A_1\to H_U\otimes_{B_s} \overline A_1$ define a morphism
$(\mathcal D_{\overline A_0},\mathcal E_{\overline A_0})\to (\mathcal D_{\overline A_1},\mathcal E_{\overline A_1})$. Thus, we have
constructed a functor $\mathcal T:YD_{brc}(\mathfrak G)\to Tens(UCorep(\mathfrak G))$.

Let now $[(\mathcal F,\eta)]:(\mathcal C_0,\mathcal E_0)\to(\mathcal C_1,\mathcal E_1)$ be a morphism in $Tens(UCorep(\mathfrak G))$,
and let $\overline A_0$ and $\overline A_1$ be the corresponding braided-commutative YD $\mathfrak G-C^*$-algebras - see Theorem
\ref{tensYD}. It follows from the construction of the $*$-algebras $A_0$ and $A_1$ that the maps $(\mathcal E_0(U^x)
\otimes\overline H^x)\to (\mathcal E_1(U^x)\otimes\overline H^x)$ given by $(X\otimes\overline\xi)\mapsto(\mathcal F(X)
\otimes\overline\xi)$ define a unital $*$-homomorphism $A_0\to A_1$ that respects their $B$-comodule and $B$-module structures.
It then extends to a homomorphism of unital braided-commutative YD $\mathfrak G-C^*$-algebras $f:\overline A_0\to \overline A_1$
which depends only on the equivalence class of $(\mathcal F,\eta)]$. Thus, we have constructed a functor $\mathcal S:Tens
(UCorep(\mathfrak G))\to YD_{brc}(\mathfrak G)$.

The homomorphism $f:\overline A_0\to \overline A_1$ is injective (resp., surjective) if and only if the maps
$Hom_{\mathcal C_0}({\bf 1},\mathcal E_0(U^x))\to Hom_{\mathcal C_1}({\bf 1},\mathcal E_1(U^x))$ are injective
(resp., surjective),
for all $x\in\Omega$. But thanks to the equalities of the type $Hom_{\mathcal C}({\bf 1},V\otimes\overline U)=Hom_{\mathcal C}(U, V)$,
this holds if and only if, for all $U,V\in UCorep(\mathfrak G)$, the maps $Hom_{\mathcal C_0}(\mathcal E_0(U),\mathcal E_0(V))\to Hom_{\mathcal C_1}(\mathcal E_1(U),
\mathcal E_1(V))$ are injective (resp., surjective). Since $\mathcal C_i$ is generated by $\mathcal E_i(UCorep(\mathfrak G))$ for $i\in \{0,1\}$, it follows
that $f$ is injective (resp., surjective) if and only if the functor $\mathcal F$ is faithful (resp., full).

%%%%%%%%%%%%%%%%%%%%%%%%%%%%%%%%%%%%%%%%%%%%%%%%%%%%%%%%%%%%%%%%%%%%%%%%%%%%%%%%%%%%%%%%%%%%%%%%%%%%%%%%%%%%%%%%%%

{\bf d) Equivalence of categories.} In order to show that the above functors $\mathcal T$ and $\mathcal S$ are inverse to each
other up to an isomorphism, let us start with a pair $(\mathcal C,\mathcal E)$ as above, the corresponding braided-commutative YD
$\mathfrak G-C^*$-algebra $\overline A$ and describe explicitly the image of any morphism $T\in Hom_{\mathcal C}(\mathcal E(U),
\mathcal E(V))$, where $U,V\in UCorep(\mathfrak G)$, under the unitary equivalence $\mathcal F:\mathcal C\to\mathcal D_{\overline A}$
as left $UCorep(\mathfrak G)$-module $C^*$-categories given by Theorem \ref{main}.

%there is a unitary equivalence $\mathcal F:\mathcal C\to\mathcal D_A$ as left $UCorep(\mathfrak G)$-module $C^*$-categories.
%First, we want to describe the image of any morphism $T\in Hom_{\mathcal C}(\mathcal E(U),\mathcal E(V))$,

In particular, $\mathcal F$ maps a morphism $T\in Hom_{\mathcal C}({\bf 1}_{\mathcal C},\mathcal E(V))=F(V)$ to the morphism
$\mathcal F(T): H^\varepsilon\otimes_{B_s}A_\varepsilon\to H_V\otimes_{B_s} A_V$ sending $1_B$ to $\underset{j}\Sigma [\zeta_j\otimes_{B_s}p_A(T\otimes\overline\zeta_j)]$, where $\{\zeta_j\}\in H_V$ is an orthonormal basis (see the proof of
\cite{VV1},Theorem 6.3). Now write any $T\in Hom_{\mathcal C}(\mathcal E(U),\mathcal E(V))$ as $T=(\mathcal E(\overline R^*_U)
\otimes id)(id\otimes S)$, where $S=(id\otimes T)\mathcal E(R_U)\in Hom_{\mathcal C}({\bf 1}_{\mathcal C},\mathcal E(\overline
U\otimes V))$. Choose an orthonormal basis $\{\xi_i\}$ in $H_U$ (as in the proof of Lemma \ref{2.8}, it is convenient to choose
$\{\xi_i\}$ such that $\hat G\cdot\xi_i=\lambda_i(U)\xi_i$ for all $i$). Then the image of the morphism $S$ is:
$$
1\mapsto\underset{i,j}\Sigma [\overline\xi_i\otimes_{B_s}\zeta_j\otimes_{B_s}p_A(S\otimes\overline{(\overline\xi_i
\otimes_{B_s}\zeta_j)})].
$$
It follows that $T=(\mathcal E(\overline R^*_U)\otimes id)(id\otimes S)$ is mapped into
$$
\underset{i,j}\Sigma [\mathcal E(\overline R^*_U) (\cdot\otimes_{B_s}\overline\xi_i)\zeta_j\otimes_{B_s}p_A(S\otimes\overline{(\overline\xi_i\otimes_{B_s}\zeta_j)})].
$$
Using the second of formulas (\ref{rigid}), we conclude that the image of $T$ is:
\begin{equation}\label{imageT}
\underset{i,j}\Sigma \theta_{\zeta_j,\xi_i}\otimes_{B_s}p_A((id\otimes T)\mathcal E(R_U)
\otimes\overline{(\hat G^{-1/2}\cdot\overline\xi_i\otimes_{B_s}\zeta_j)},
\end{equation}
where $\theta_{\zeta_j,\xi_i}\in B(H_U,H_V)$ is defined by $\theta_{\zeta_j,\xi_i}(\eta)=<\eta,\xi_i>\zeta_j$ for all $\eta\in H_U$.

In order to show that $\mathcal F$ is a strict tensor functor on $\mathcal E(UCorep(\mathfrak G))$ and hence on $\mathcal C$, we have
to show that $\mathcal F(S\otimes T)=\mathcal F(S)\otimes \mathcal F(T)$ on morphisms in $\mathcal E(UCorep(\mathfrak G))$. Since
$\mathcal F$ is an equivalence of left $UCorep(\mathfrak G)$-module categories, we already know that $\mathcal F(id\otimes T)=
id\otimes \mathcal F(T)$, so it remains to show that $\mathcal F(S\otimes id)=\mathcal F(S)\otimes id$.

If $S:\mathcal E(U)\to \mathcal E(V)$ and $\{\eta_k\}$ is an orthonormal basis in $H_W\ (W\in UCorep(\mathfrak G))$, then according
to (\ref{imageT}) $\mathcal F(S\otimes id_W)$ equals
$$
\underset{i,j,k,l}\Sigma \theta_{\zeta_j\otimes_{B_s}\eta_k,\xi_i\otimes_{B_s}\eta_l}\otimes_{B_s}p_A((id\otimes (S\otimes id))
\mathcal E(R_{U\otimes W})\otimes
$$
$$
\otimes\ \overline{\hat G^{-1/2}\cdot\overline{(\xi_i\otimes_{B_s}\eta_l)}\otimes_{B_s}\zeta_j\otimes_{B_s}\eta_k}).
$$
As in the proof of Lemma \ref{2.6}, $R_{U\otimes W}$ coincides, modulo the equivalence  $\overline{U\otimes W}\equiv
\overline W\otimes\overline U$, with $(id\otimes\overline R_U\otimes id)\overline R_W)$, so the above expression equals
$$
\underset{i,j,k,l}\Sigma \theta_{\zeta_j\otimes_{B_s}\eta_k,\xi_i\otimes_{B_s}\eta_l}\otimes_{B_s}p_A((id\otimes(id\otimes S)
\mathcal E(R_U)\otimes id)\mathcal E(R_{W})\otimes
$$
$$
\otimes\overline{\hat G^{-1/2}\cdot\overline\eta_l\otimes_{B_s}\hat G^{-1/2}\cdot\overline\xi_i\otimes_{B_s}\zeta_j\otimes_{B_s}
\eta_k}).
$$
The operators $\theta_{\eta_k,\eta_l}$ are the matrix units $m_{k,l}$ in $B(H_W)$. Recalling the definition of $\triangleleft$,
we can rewrite the above expression as
$$
\underset{i,j,k,l}\Sigma \theta_{\zeta_j,\xi_i}\otimes_{B_s}m_{k,l}\otimes_{B_s}[p_A((id\otimes S)\mathcal E(R_U)
\otimes\overline{(\hat G^{-1/2}\cdot\overline\xi_i\otimes_{B_s}\zeta_j)})\triangleleft W_{k,l}]
$$
On the other hand, $\mathcal F(S):H_U\otimes_{B_s} A\to H_V\otimes_{B_s} A$ can be presented as $\mathcal F(S)=\underset{i}
\Sigma s_i\otimes_{B_s} a_i$, where $s_i\in B(H_U,H_V), a_i\in A$, with the action $\mathcal F(S)(\xi\otimes_{B_s} a)=
\underset{i}\Sigma s_i(\xi)\otimes_{B_s} a_i a$, for all $\xi\in H_U, a\in A$. Considering $\mathcal F(S)\otimes id$ as a
morphism from $H_{U\otimes W}\otimes_{B_s} A$ to $H_{V\otimes W}\otimes_{B_s} A$, we have for all $\zeta\in H_U, \eta\in H_W$:
$$
(\mathcal F(S)\otimes id)(\zeta\otimes_{B_s}\eta\otimes_{B_s}1_{\overline A})=T_{V,W}(\underset{i}\Sigma s_i(\zeta
)\otimes_{B_s}\eta\otimes_{B_s} a_i)=
$$
$$
=T_{V,W}(\underset{i}\Sigma s_i(\zeta)\otimes_{B_s}a_i)\otimes_{\overline A}(\eta\otimes_{B_s}1_{\overline A})=
$$
$$
\underset{i}\Sigma s_i(\zeta)\otimes_{B_s}\underset{k,l}\Sigma m_{k,l}\eta\otimes_{B_s}(a_i\triangleleft W_{k,l}).
$$
Hence, the actions of  $\mathcal F(S\otimes id)$ and $\mathcal F(S)\otimes id$ on generating vectors $\zeta\otimes_{B_s}
\eta\otimes_{B_s}1_{\overline A}$ coincide.
\begin{remark} \label{comput}
Similar calculation and the fact that $1_{\overline A}\triangleleft U_{k,l}=\delta_{k,l}1_{\overline A}$ give
$$
(id\otimes\mathcal F(S))(\eta\otimes_{B_s}\zeta\otimes_{B_s}1_{\overline A})=\eta\otimes_{B_s}\underset{i}\Sigma s_i(\zeta)
\otimes_{B_s}a_i.
$$
\end{remark}
Conversely, consider a unital braided-commutative YD $\mathfrak G-C^*$-algebra $\overline A$ and the corresponding pair
$(\mathcal C_{\overline A},\mathcal E_{\overline A})$, and let $\overline A_{\mathcal C}$ be the braided-commutative YD
$\mathfrak G-C^*$-algebra constructed from this pair. By Theorem \ref{main}, there is an isomorphism $\lambda:\overline
A_{\mathcal C}\to \overline A$ intertwining the coactions of $\mathfrak G$ and defined by $\lambda(p_{\overline A}(T\otimes\overline\zeta))=(\overline\zeta\otimes id)T$, for all $\zeta\in H_V, T\in \mathcal C_{\overline A}({\bf 1},V)
\subset\mathcal L (B_s,H_V)\otimes{\overline A} =H_V\otimes {\overline A} $.
So it only remains to show that $\lambda$ is a right $B$-module map.

As above, fix $U,V\in UCorep(\mathfrak G)$ and orthonormal bases $\{\xi_i\}\in H_U$ and $\{\zeta_k\}\in H_V$, and let
$U_{kl}$ be the matrix coefficients of $U$. Take $T=\underset{k}\Sigma(\zeta_k\otimes_{B_s} a_k)\in H_V\otimes_{B_s} A$,
then $\lambda(p_{\overline A}(T\otimes\overline\zeta_{k_0}))=a_{k_0}$, and check that $\lambda(p_{\overline A}
(T\otimes\overline\zeta_{k_0})\triangleleft U_{i_0,j_0})=a_{k_0}\triangleleft U_{i_0,j_0}$. By (\ref{module}) we have
$$
p_{\overline A}(T\otimes\overline\zeta_{k_0})\triangleleft U_{i_0,j_0}=
%p_{\overline A}(T\otimes\overline\zeta_{k_0})\triangleleft p_B(\xi_{i_0}\otimes_{B_s}\xi_{j_0})=
(id\otimes T\otimes id)F(R_U)\otimes\overline{(\hat G^{-1/2}\overline\xi_{i_0}
\otimes_{B_s}\zeta_{k_0}\otimes_{B_s}\xi_{j_0})}.
$$
%where $\{\zeta_k\}$ is an orthonormal base in $H_V, a_k\in A$, and the formula (\ref{module}) for the action of $B$ on $A$.
In order to compute the image of this element under $\lambda$, we need an explicit formula for $(id\otimes T\otimes id)F(R_U)
\ :\ {\bf 1}\to H_{\overline U\otimes V\otimes U}\otimes_{B_s} A$. Remark \ref{comput} and the computation before it show that the
element $id_{\overline U}\otimes T\otimes id_U: H_{\overline U}\otimes_{B_s} H_U\to H_{\overline U}\otimes_{B_s} H_V \otimes_{B_s}
H_U$ equals $1\otimes_{B_s}\underset{i,j,k}\Sigma\zeta_k\otimes_{B_s}m_{i,j}\otimes_{B_s}(a_k\triangleleft U_{j,i})$. Then
$$
(id\otimes T\otimes id)F(R_U)=\underset{i,j,k}\Sigma\hat G^{1/2}\cdot\overline\xi_j\otimes_{B_s}\zeta_k\otimes_{B_s}\xi_i
\otimes_{B_s}(a_k\triangleleft U_{j,i}).
$$
Therefore, $p_{\overline A}(T\otimes\overline\zeta_{k_0})\triangleleft U_{i_0,j_0}$ equals
$$
p_{\overline A}(\underset{i,j,k}\Sigma\hat G^{1/2}\cdot\overline\xi_j\otimes_{B_s}\zeta_k\otimes_{B_s}\xi_i\otimes_{B_s}
(a_k\triangleleft U_{j,i})\otimes\overline{(\hat G^{-1/2}\overline\xi_{i_0}\otimes_{B_s}\zeta_{k_0}\otimes_{B_s}\xi_{j_0})})
$$
Applying $\lambda$, we get the required equality $\lambda(p_{\overline A}(T\otimes\overline\zeta_{k_0})=a_{k_O}\triangleleft
U_{i_0,j_0}$.
As the algebra $A$ is spanned by elements $a_{k_0}$ for various $V\in UCorep(\mathfrak G)$, it follows that $\lambda$ is a right
$B$-module map. This completes the proof of Theorem \ref{catdual}.

%\begin{remark} \label{ergodic}
%It follows from \cite{VV1}, Remark 4.11, that if $(A,\mathfrak a)$ is a braided-commutative YD $C^*$-algebra, then the
%coaction $\mathfrak a$ is necessarily ergodic.
%\end{remark}
\end{subsection}
\end{section}

%%%%%%%%%%%%%%%%%%%%%%%%%%%%%%%%%%%%%%%%%%%%%%%%%%%%%%%%%%%%%%%%%%%%%%%%%%%%%%%%%%%%%%%%%%%%%%

\begin{section} {Quotient type and invariant coideals}

\begin{subsection} {Quotient type coideals} The notion of a quotient type coideal of a WHA $\mathfrak G$ is closely
related to the notion of a quantum subgroupoid which is just another WHA $\mathfrak H$ equipped with an epimorphism
$\pi:\mathfrak G\to \mathfrak H$. We start with basic definitions and results.

\begin{definition} \label{morphism} A morphism of two WHAs, $\mathfrak G=(B,\Delta^B,S^B,\varepsilon^B)$
and $\mathfrak H=(C,\Delta^C,S^C,\varepsilon^C)$, is a unital morphism $\pi:B\to C$ of their $C^*$-algebras
such that $\Delta^C\circ\pi=(\pi\otimes\pi)\Delta^B$, $S^C\circ\pi=\pi\circ S^B$ and $\varepsilon^C\circ\pi=\varepsilon^B$.
\end{definition}

\begin{remark} \label{pi} 1. One checks that this definition implies: $\pi\circ\varepsilon^B_t= \varepsilon^C_t\circ\pi$
and $\pi\circ\varepsilon^B_s=\varepsilon^C_s\circ\pi$,
%so $\pi(B_t)\subset C_t$ and $\pi(B_s)\subset C_s$.
so $\pi(B_t)= C_t$ and $\pi(B_s)= C_s$.

2. If $\mathfrak G$ and $\mathfrak H$ are usual Hopf $C^*$-algebras, Definition \ref{morphism}
%gives $\varepsilon^B(\pi(1_B))=\varepsilon^C(\pi(b)1_C$, for all $b\in B$, so it
coincides with the usual definition of a morphism of Hopf $C^*$-algebras.
%if and only if either $\pi(1_B)=1_C$ or $\varepsilon^B=\varepsilon^C\circ\pi$.
\end{remark}

\begin{lemma} \label{tfunctor}
If $\pi:\mathfrak G\to \mathfrak H$ is surjective, the map $\mathcal E_{\pi}:U\mapsto (id\otimes\pi)U$
is a unitary tensor functor $UCorep(\mathfrak G)\to UCorep(\mathfrak H)$. Moreover, $(id\otimes\pi)U\in
B(H^\pi_U)\otimes C$, where $(H^\pi_U)^\perp=\{\zeta\in H_U|(id\otimes\pi)U(\zeta\otimes 1)=0\}$.
\end{lemma}
\begin{dm}
Considering $H_U$ as a left $\hat B$-module in the following way:
$$
\hat b\cdot\zeta=<U^{(2)},\hat b>U^{(1)}\zeta\quad(\forall \hat b\in \hat B,\zeta\in H_U),
$$
and the $*$-algebra inclusion $\pi_*:\hat C\to\hat B$ dual to $\pi:B\to C$, one has:
$$
<\pi(U_{\eta,\zeta}),\hat c>=<(id\otimes\pi)U(\zeta\otimes 1_B),\eta\otimes\hat c>=<U^{(1)}\zeta,\eta><U^{(2)},\pi_*(\hat c)>=
$$
$$
=<\pi_*(\hat c)\zeta,\eta>=<\zeta,\pi_*(\hat c)^*\eta>\quad (\forall \eta,\zeta\in H_U,\hat c\in\hat C),
$$
from where $H^\pi_U=\pi_*(\hat C)H_U$. In particular, $\pi(U_{\eta,\zeta})=0$ for all $\zeta\in H^\pi_U, \eta\in (H^\pi_U)^\perp$
which gives the result.\hfill\end{dm}

\begin{corollary} \label{ctfunctor} The functor $\mathcal E_\pi$ transforms $H_U$ into $H^\pi_U$ and
intertwiners $H_U\to H_V$ into intertwiners $H^\pi_U\to H^\pi_V$ for all $U,V\in UCorep(\mathfrak G)$.
\end{corollary}

\begin{lemma} \label{matrcoef} Let $\hat l\in\hat C$. The matrix coefficient $\pi(U_{\pi_*(\hat l)\zeta,\eta})
\in C_s$ for all $U\in UCorep(\mathfrak G),\eta,\zeta\in H_U$ if and only if $\hat l$ is a left integral.
\end{lemma}

\begin{dm} Combining the above mentioned relations, first one has:
$$
\pi(U_{\eta,\pi_*(\hat l)\zeta})=<U^{(1)}\pi_*(\hat l)\zeta,\eta>\pi(U^{(2)})=
$$
$$
=<U^{(1)}U^{'(1)}\zeta,\eta>\pi(U^{(2)})<U^{'(2)},\pi_*(\hat l)>=
$$
$$
<U^{(1)}\zeta,\eta>(id\otimes\hat l)\Delta_C(\pi(U^{(2)})).
$$
Since $C$ is spanned by the $\pi$-images of matrix elements of all $U\in UCorep(\mathfrak G)$, this element is
in $C_t$ if and only if $(id\otimes\hat l)\Delta_C=(\varepsilon^C_t\otimes\hat l)\Delta_C$, i.e., if and only if
$\hat l$ is a left integral. Finally, $\pi(U_{\pi_*(\hat l)\zeta,\eta})=S_C(\pi(U_{\eta,\pi_*(\hat l)\zeta}))^*$.
\hfill\end{dm}
\vskip 0.5cm
Consider now a surjective $*$-homomorphism $\pi : \mathfrak G\twoheadrightarrow \mathfrak H$ of
WHAs $\mathfrak G=(B,\Delta_B,S_B,\varepsilon_B)$ and $\mathfrak H=(C,\Delta_C,S_C,\varepsilon_C)$.
%such that $\varepsilon_C\circ\pi=\varepsilon_B$.
Then the map $\mathfrak a=(\pi\otimes id_B)\Delta$ is a {\em left} coaction of $\mathfrak H$ on $B$.
\begin{definition} \label{qtype}
The fixed point unital $*$-subalebra $I(\mathfrak H\backslash\mathfrak G)$ of $B$ with respect to the coaction
$(\pi\otimes id_B)\Delta$ of $\mathfrak H$ is called {\it a quotient type coideal subalgebra} (briefly, quotient type coideal)
of $B$. Equivalently,
%one can check that the subset $I(\mathfrak H\backslash\mathfrak G)\subset B$ defined by:
$$
I(\mathfrak H\backslash\mathfrak G) = \{b \in B | {} (\pi\otimes id)\Delta_B(b) = (\pi\otimes id)((1\otimes b)\Delta_B(1)) \}
$$
%is a right coideal subalgebra (called {\it a quotient type coideal subalgebra}) of $B$.
\end{definition}
Obviously, $\pi(I)\subset C_s$ and $B_s\subset I$, so $\pi(I)= C_s$.
%Note also that $I$ is exactly the fixed point subalgebra of $B$ with respect to the coaction $(\pi\otimes id)\Delta$ of $\mathfrak H$.

Clearly, the smallest quotient type coideal of $B$ is $B_s$. It corresponds to $\mathfrak H=\mathfrak G$, $\pi=id$.
Since $I(\mathfrak H\backslash\mathfrak G)$ is the fixed point $*$-subalgebra with respect to the coaction $(\pi\otimes id)\Delta$,
it is included into $B'_t\ (=\alpha(B_t)')$.

%Lemma \ref{quinv} shows that any quotient type coideal is the subset of $B'_t$. More exactly, we have
\begin{lemma} \label{minimal}
$B_t'$ is the greatest quotient type coideal.
\end{lemma}
\begin{dm}
Let $\hat{\mathfrak G}$ be the dual of a WHA $\mathfrak G=(B,\Delta,S,\varepsilon)$. Let
$\hat{\mathfrak G}_{min}=\hat B_t\hat B_s$ be the minimal WHA contained in $\hat{\mathfrak G}$ and $i_{min}:
\hat{\mathfrak G}_{min}\to\hat{\mathfrak G}$ the corresponding inclusion of WHAs (see \cite{BNSz}, \cite{Nik}).
Then the adjoint map $\pi_{min}:\mathfrak G\to\hat{\mathfrak G}^*_{min}$ given by
$<i_{min}(\hat B_{min}),B>=<\hat B_{min},\pi_{min}(B)>$, is an epimorphism of WHAs. The corresponding quotient type
coideal $I_{max}$ is the set of such $b\in B$ that
$$
<\hat z\otimes\hat b,(\pi_{min}\otimes id)\Delta(b)>=<\hat z\otimes\hat b,(\pi_{min}\otimes id)(1_B\otimes b)\Delta(_B1)>,
$$
which is equivalent to
\begin{equation} \label{min}
<\hat b,b\leftharpoonup \hat z>=<\hat b,b(1_B\leftharpoonup \hat z)>,\quad\text{for\ all}\quad \hat b\in\hat B,\ \hat z\in
{\hat B}_{min},
\end{equation}
where, by definition,
$c\leftharpoonup \hat z:=c_2<\hat z,c_1>$, for any $c\in B$. As ${\hat B}_{min}={\hat B}_s{\hat B}_t$, it suffices to consider
$\hat z=\hat u\hat v$, where $\hat u\in {\hat B}_t$ and $\hat v\in {\hat B}_s$. So, $b\in I_{max}$ if and only if $b\in B$ and
$$
b\leftharpoonup ({\hat u}{\hat v})=b(1_B\leftharpoonup ({\hat u}{\hat v}))b\quad \text{or}\quad
(b\leftharpoonup{\hat u})\leftharpoonup \hat v=b((1_B\leftharpoonup{\hat u})\leftharpoonup \hat v).
$$
By \cite{BNSz}, (2.21a), one can rewrite the last equality as $(ub)\leftharpoonup(\hat v)=b(u\leftharpoonup\hat v)$, where we
denoted $1_B\leftharpoonup \hat u\in B_t$ by $u$. Now, by \cite{BNSz}, (2.20b), this can be rewritten as $(ub)v=b(uv)$, where we
denoted $1_B\leftharpoonup \hat v\in B_s$ by $v$. But  this is true exactly for $b\in B'_t$, and we are done.
\hfill\end{dm}

%\begin{remark} \label{elem} Duals of minimal WHA's were considered in \cite{NV5}.
%\end{remark}

\begin{lemma} \label{mc} Let $\zeta,\eta\in H_U$, then $U_{\zeta,\eta}\in I$ if and only if $\pi(U_{\zeta,\theta})\in C_s$,
for all $\theta\in H_U$.
\end{lemma}
\begin{dm} If $\{\theta_k\}$ is an orthonormal basis is $H_U$ and $U_{\zeta,\eta}\in I$, then $\Delta(U_{\zeta,\eta})=
\underset{k}\Sigma (U_{\zeta,\theta_k}\otimes U_{\theta_k,\eta})\in I\otimes B$, so all $U_{\zeta,\theta_k}\in I$ and
$\pi(U_{\zeta,\theta_k})\in C_s$ for all $k$. So, $\pi(U_{\zeta,\theta})\in C_s$, for all $\theta\in H_U$.

Conversely, if $\pi(U_{\zeta,\theta})\in C_s$, then
$$
(\pi\otimes id)\Delta_B(U_{\zeta,\theta})=(\varepsilon^C_s\otimes id)(\pi\otimes id)\Delta_B(U_{\zeta,\theta})=
$$
$$
=(\pi\otimes id)(\varepsilon^B_s\otimes id)\Delta_B(U_{\zeta,\theta})=(\pi\otimes id)(1\otimes U_{\zeta,\theta})\Delta_B(1).
$$
\end{dm}

The following lemma generalizes \cite{VV1}, Example 6.7 and describes module categories associated
with quotient type coideals.

\begin{lemma} \label{genregular} If $(\mathfrak H,\pi)$ is a quantum subgroupoid of $\mathfrak G$ and
$\Lambda$ is the set of left integrals of $\hat{\mathfrak H}$, then $I=
I(\mathfrak H\backslash\mathfrak G)$ admits the decomposition
$$
I(\mathfrak H\backslash\mathfrak G)=\underset{x\in\Omega}\oplus\ \pi_*(\Lambda)H^x\otimes H^x
$$
and the corresponding $UCorep(\mathfrak G)$-module $C^*$-category $\mathcal M$ is equivalent to
$UCorep(\mathfrak H)$ viewed as a $UCorep(\mathfrak G)$-module $C^*$-category via the functor $\mathcal E_\pi$.
\end{lemma}
\begin{dm} It suffices to prove that $I$ is equivariantly isomorphic to the $\mathfrak G-C^*$-algebra $A$
corresponding to the couple $(UCorep(\mathfrak H),\bf 1)$. Following the categorical duality, we first construct
an algebra $\tilde A$ of the form (\ref{recalgebra2}), where $F(U)=Hom_{UCorep(\mathfrak H)}({\bf 1},\mathcal
E_\pi(H_U))=Hom_{URep(\hat{\mathfrak H})}(C_s,\pi_*(\hat C)(H_U))$, where $C_s$ is a left $\hat C$-module via
$\hat c\cdot z:=\hat c\rightharpoonup z\ (\hat c\in \hat C,z\in B_s)$. For any such morphism $f$ the vector $f(1)$
is cyclic for $Im(f)$, so in fact we have to describe $Hom_{Rep(\hat{\mathfrak H})}(C_s,\hat C)$. But \cite{BNSz},
Lemma 3.3 shows that this is exactly the set $\Lambda$ of left integrals in $\hat C$. Thus, we can identify
$Hom_{URep(\hat{\mathfrak H})}(C_s,\pi_*(\hat C)(H_U))$ with the subspace $\pi_*(\Lambda)H_U\subset H_U$.

So $\tilde A$ is a subalgebra of the algebra $\tilde B=\underset{U}\bigoplus (H_U\otimes\overline H_U)$, and the map
$p_B:\tilde B\to B$ sending $\zeta\otimes\overline\eta\in H_U\otimes\overline H_U$ onto the matrix coefficient $U_{\zeta,\eta}$
induces an $\mathfrak G$-equivariant isomorphism of $A$ onto the coideal $I=Vec\{U_{\zeta,\eta}|\zeta\in \pi_*(\Lambda)H_U,
\eta\in H_U, U\in UCorep(\mathfrak G)\}$. Finally, Corollary \ref{matrcoef} and Lemma \ref{mc} show that
$I=I(\mathfrak H\backslash\mathfrak G)$.
\hfill\end{dm}

\end{subsection}

\begin{subsection} {Invariant coideals}

Consider the {\bf right adjoint action} of $B$ on itself defined by
\begin{equation} \label{adjoint}
b\triangleleft x:=S(x_{(1)}) b x_{(2)},\ \text{for\ all}\ x,b\in B.
\end{equation}
It follows from \cite{NTV}, Lemma 2.2 that
%As $(S\otimes id)\Delta(1)$ is a separability element for $B_t$ - see \cite{NV}, Prop. 3.2.4,
the map $P_\triangleleft :b\mapsto b\triangleleft 1_B$ is a projection from $B$ onto $B'_t$,
from where $B\triangleleft B =B_t'$.

\begin{definition} \label{invar}
A right coideal $I$ is called invariant if $I\triangleleft B= I$.
\end{definition}

\begin{remark} \label{reminv}
1. $I$ is invariant if and only if $I\subset B'_t$ and $I\triangleleft B\subset I$.
Indeed, Definition \ref{invar} implies $I=I\triangleleft B\subset B\triangleleft B =B_t'$ and
$I\triangleleft B=I\subset I$. Conversely, if $I\subset B'_t$ and $I\triangleleft B\subset I$,
then $I=P_\triangleleft(I)=I\triangleleft 1_B\subset I\triangleleft B\subset I$.

2. Corollary: $B$ is invariant if and only if $\mathfrak G$ is a $*$-Hopf algebra.

3. One can check that $B'_t$ is the greatest invariant coideal.
%and $B_s$ is the smallest {\bf (?)}
%As any invariant weak coideal contains $1_B\in B_s$, it is a coideal.
\end{remark}

It is known \cite{NV5} that all coideals of $B$ form a lattice $l(B)$ with minimal element $B_s$
and maximal element $B$ under the usual operations:
$I_1 \wedge I_2 = I_1 \cap I_2,\ I_1 \vee I_2 = (I_1 \cup I_2)''$.

\begin{lemma} \label{lattice }
Invariant coideals form its sublattice $invl(B)$ with minimal (resp., maximal) element $B_s$ (resp., $B'_t$).
\end{lemma}
\begin{dm}
If $I,I'$ are two invariant  coideals, then for any natural number $k\geq 2$ and all $j_1,.. ,j_k
\in I \cup I', b \in B$ one has:
$$
j_1..j_k \triangleleft b = ( j_1...j_{k-1}\triangleleft b_{(1)})(j_k\triangleleft b_{(2)}).
$$
So by obvious iteration $j_1..j_k \triangleleft b \in I\vee I'$, but $I,I'$ are $*$-invariant so $I \vee I'$
is spanned by sums of type $j_1..j_k$. Moreover, the map: $j \mapsto j \triangleleft b$ is linear which gives that
$(I \vee I')\triangleleft B\subset (I \vee I')$. But $(I \vee I')\subset B'_t$, so by Remark \ref{reminv}, 1
$I \vee I'$ is invariant. Also, $(I \cap I')\triangleleft B\subset (I \cap I')$ and $(I \cap I')\subset B'_t$, so
$I \cap I'$ is invariant and the result follows.\hfill
\end{dm}

\begin{lemma} \label{invYD}
Any invariant coideal $I$ belongs to the category $YD_{brc}(\mathfrak G)$.

%is a braided-commutative YD $\mathfrak G$-$C^*$-algebra.
%{\bf (Only if $B_s$ is abelian)}
\end{lemma}
\begin{dm}
All the relations (\ref{triangle}) are obvious.
%apart from the last one. But since
%$\Delta(1_I)\in I\otimes B_t$, we have, using \cite{NV}, Prop. 2.2.1, that $1_I\in B_t$, so
%$$
%1_I\triangleleft b=S(b_{(1)})1_I b_{(2)}=S(b_{(1)}S^{-1}(1_I))b_{(2)}=1_I\varepsilon_s(b),\ \forall b\in B.
%$$
%The needed equality follows from this and the next calculation:
%$$
%1_I\triangleleft\varepsilon_s(b)=S(1_{(1)})1_I\varepsilon_s(b)1_{(2)}=1_I\varepsilon_s(b),\ \forall b\in B.
%$$
Let us check the Yetter-Drinfel'd and the braided commutativity relations:
$$
\Delta(S(b_{(1)})ab_{(2)})=S(b_{(2)})a_{(1)}b_{(3)}\otimes S(b_{(1)})a_{(2)}b_{(4)}=
(a_{(1)}\triangleleft b_{(2)})\otimes S(b_{(1)})a_{(2)}b_{(3)},
$$
where $a\in I, b\in B$. Finally, using the fact that $I\in B'_t$:
$$
b_{(1)}(a\triangleleft b_{(2)})=b_{(1)}S(b_{(2)})ab_{(3)}=\varepsilon_t(b_{(1)})ab_{(2)}=
a\varepsilon_t(b_{(1)})b_{(2)}=ab,\ \forall a,b\in I.
$$
\hfill\end{dm}

Now discuss the relationship between quotient type and invariant coideals.

\begin{lemma} \label{quinv} Any quotient type coideal $I$ is invariant.
\end{lemma}
\begin{dm}
Let us show that $(I\triangleleft x)\in I$ for all $x\in B$.
Indeed, using \cite{NV}, Proposition 2.2.1, we have for all $b\in I$:
\begin{align*}
(\pi\otimes& id)\Delta(S(x_{(1)})b x_{(2)}))=\\
&=(\pi(S(x_{(2)})\otimes S(x_{(1)}))((\pi\otimes id)((1\otimes b)
\Delta(1))(\pi(x_{(3)})\otimes x_{(4)})=
\\
&=(\pi\otimes id)(S(x_{(2)})1_{(1)} x_{(3)}\otimes S(x_{(1)})b 1_{(2)}x_{(4)})=
\\
&=(id\otimes\pi)(S(x_{(2)})x_{(3)}\otimes S(x_{(1)}) b S(x_{(4)})))=\\
&=(\pi\otimes id)(\varepsilon_s(x_{(2)})\otimes S(x_{(1)})b x_{(3)})
=(\pi\otimes id)(1_{(1)}\otimes S(x_{(1)}) b x_{(2)} 1_{(2)}).
\end{align*}
\hfill\end{dm}

The inverse statement - Theorem \ref{inv-qt} is proved as follows:
%Similarly to the compact quantum group case \cite{NY1}, Theorem 3.1, we have

%\begin{theorem} \label{inv-qt}
%For any invariant weak coideal $I$
%is covariantly isomorphic to a unique quotient type coideal, i.e.,
%there exists a unique (up to isomorphism) quantum subgroupoid $\mathfrak H$ such that $I\cong I(\mathfrak H\backslash \mathfrak G)$
%as $\mathfrak G$-$C^*$-algebras.
%\end{theorem}

\begin{dm} Lemma \ref{invYD} shows that $I$ is a braided-commutative YD $\mathfrak G$-$C^*$-algebra. Then Theorem
\ref{2.4} shows that the corresponding $UCorep(\mathfrak G)$-module category $\mathcal D_I$ is a $C^*$-multitensor
category with tensor product $\otimes_I$ and trivial associativities equipped with a unitary tensor functor
$$\mathcal E_I: UCorep(\mathfrak G) \to \mathcal D_I \ \ , \ \  U \mapsto   H_U\otimes_{B_s} I.$$

Let us equip now the category $\mathcal D_I$ with the tensor functor $\mathcal F$ to $Corr_f(B_s)$ sending $H_U
\otimes_{B_s} I$ to $H_U$. Then the reconstruction theorems for WHA's and their morphisms allow to construct a WHA
$\mathfrak H_I$ such that $UCorep(\mathfrak H_I)\cong\mathcal C_I$ together with an epimorphism $\pi:\mathfrak G\to
\mathfrak H_I$. In its turn, this allows to construct the quotient type coideal $J=I(\mathfrak H_I\backslash
\mathfrak G)$. Now, Lemma \ref{genregular} shows that $\mathcal D_{J}\cong UCorep(\mathfrak H_I)$, therefore,
$\mathcal D_{J}\cong \mathcal D_I$. But due to the categorical duality this implies a covariant isomorphism $J\cong I$.

In order to prove the uniqueness of $\mathfrak H$, we prove that $J$ determines $Ker(\pi)$. Indeed, as $J=Vec\{U_{\eta,
\zeta}\}$, where $U\in UCorep(\mathfrak G)\}, \eta\in H_U$ and $\zeta\in\pi_*(\hat L)H_U\subset H_U$,
the last subspace is determined by $J$. This means that for any $U\in UCorep(\mathfrak G$, the space $Hom_{UCorep
(\mathfrak H)}({\bf 1},\mathcal (H^\pi_U))$ is determined by $J$. The duality morphisms $Hom_{UCorep(\mathfrak G)}
({\bf 1},U\otimes \overline V)=Hom_{UCorep(\mathfrak G)}(V,U)$ $(U,V\in UCorep(\mathfrak G))$ show that the same is true
for all subspaces \newline $Hom_{UCorep(\mathfrak H}(H^\pi_V),H^\pi_U)\subset B(H_V,H_U)$. Finally, an operator from
$B(H_U)$ belongs to $Ker(\pi)$ if and only if its matrix coefficients corresponding to the commutant of $End_{UCorep
(\mathfrak H}(H^\pi_U)$ in $B(H_U)$ equal to $0$.
\hfill\end{dm}

\begin{remark} \label{weak} In \cite{VV2} we introduced the notion of a {\em weak} coideal $I$ of $B$, the difference
of which from a coideal is that $1_I$ is not necessarily equal to $1_B$. One can show that if $I$ is invariant and
$1_I\in Z(B)$, then $(I,\Delta)\in YD_{brc}(\mathfrak G)$. Then the same reasoning as in the proof of Theorem \ref{inv-qt}
shows that $I$ is isomorphic, as a $\mathfrak G$-$C^*$-algebra, to a unique quotient type coideal.
\end{remark}

\end{subsection}
\end{section}

%%%%%%%%%%%%%%%%%%%%%%%%%%%%%%%%%%%%%%%%%%%%%%%%%%%%%%%%%%%%%%%%%%%%%%%%%%%%%%%%%%%%%%%%%%%%%%%%%%%%%%%%%%%%%%%%%%%%%%

\begin{section} {Example: the Tambara - Yamagami case}

%Now let's illustrate the situation when $\mathcal C$ is not only multitensoriel but tensoriel.

\begin{remark} \label{rec} In this example $\mathcal C$ is not only $C^*$-multitensor, but a rigid finite $C^*$-tensor category.
Let $\mathcal F : \mathcal C \to Corr_f(R)$ be a unitary tensor functor, where $R$ is a finite dimensional unital $C^*$-algebra.
Then it was shown in \cite{Sz} that the WHA reconstructed from the pair $(\mathcal C,\mathcal F)$ as in Theorem \ref{genreconstruction}
is biconnected: $B_t\cap B_s=\mathbb C=B_t\cap Z(B)$. Moreover, Hayashi \cite{Ha} proved that for any given $\mathcal C$, the class of
$C^*$-algebras $R$ for which such a functor $\mathcal F$ exists, contains at least $R=\mathbb C^{|\Omega|}$ (where $\Omega=
Irr(\mathcal C)$). He also constructed the corresponding particular functor $\mathcal H$. In general, this class
of $C^*$-algebras $R$ contains several elements, and the corresponding WHAs are called Morita equivalent. In particular, if
this class contains
%a {\bf fiber} functor $\mathcal F : \mathcal C \to Hilb_f$, i.e.,
$R=\mathbb C$, the corresponding WHAs are Morita equivalent to a usual $C^*$-Hopf algebra.
\end{remark}

\begin{subsection} {Reconstruction for Tambara-Yamagami categories}

The description of the Hayashi's functor for Tambara-Yamagami categories and the corresponding WHA's was originally
obtained in \cite{M}. Below we follow \cite{VV2}, 2.3 and 4.1, where one can find more details.

Given a finite abelian group $G$, a non degenerate symmetric bicharacter $\chi$ on it and a number $\tau=\pm |G|^{-1/2}$, one can
define a fusion category denoted by $\mathcal T\mathcal Y(G,\chi,\tau)$ \cite{TY}. Its set of simple objects is $\Omega = G \sqcup
\{m\}\ (m$ is a separate element), its Grothendieck ring is isomorphic to the $\mathbb Z_2$-graded fusion ring $\mathcal T\mathcal Y_G
=\mathbb Z G\oplus\mathbb Z\{m\}$ such that $g\cdot m=m\cdot g=m,\ m^2=\underset{g\in G}
\Sigma g, \ g^*=g^{-1}, \ m=m^*$.  The associativities
$a_{U,V,W}:(U\otimes V)\otimes W\to U\otimes(V\otimes W)$ are
$$
a_{g,h,k}=id_{g+h+k},\quad a_{g,h,m}=id_{m},\quad a_{m,g,h}=id_{m},
$$
$$
a_{g,m,g}=\chi(g,h)id_{m},\quad a_{g,m,m}=\underset{h\in G}\oplus id_{h},\quad a_{m,m,g}=\underset{h\in G}\oplus id_{h},
$$
$$
a_{m,g,m}=\underset{h\in G}\oplus \chi(g,h)id_{h},\quad a_{m,m,m}=(\tau\chi(g,h)^{-1} id_{m})_{g,h},
$$
where $g,h,k\in G$. The unit isomorphisms are trivial. $\mathcal T\mathcal Y(G,\chi,\tau)$ becomes a $C^*$-tensor category
when $\chi:G\times G\to T=\{z\in\mathbb C||z|=1\}$, from now on we assume that this is the case. The dual objects are: $g^*=-g$,
for all $g\in G$, and $m^*=m$. The rigidity morphisms are defined by $R_g:0\overset{id_0}\to g^*\otimes g$, $\overline{R_g}:
0\overset{id_0}\to g\otimes g^*$, $R_m=\tau |G|^{1/2}\iota$, and $\overline{R_m}=|G|^{1/2}\iota$, where $\iota:0\to
m\otimes m$ is the inclusion. Then $dim_q(g)=1$, for all $g\in G$, and $dim_q(m)=\sqrt{|G|}$.

%In\cite{M} is constructed
Using now the Hayashi's functor $\mathcal H:\mathcal T\mathcal Y(G,\chi,\tau)\to Corr_f(R)$, where $R\cong\mathbb C^{|G|+1}$
(see \cite{M}, \cite{VV2}), one can apply Theorem \ref{genreconstruction} in order to construct a biconnected regular WHA
$\mathfrak G_{\mathcal T\mathcal Y}=(B,\Delta,S,\varepsilon)$ with $UCorep(\mathfrak G_{\mathcal T\mathcal Y})\cong\mathcal T
\mathcal Y(G,\chi,\tau)$ as $C^*$-tensor categories. It happens that $\mathfrak G_{\mathcal T\mathcal Y}$ is selfdual.

Denoting $\Omega_g=\Omega:=G\sqcup\{m\}$ and $\Omega_m:=G\sqcup\overline G$, where $g\in G$ and $\overline G$ is the second
copy of $G$, one computes that $H^g\cong\mathbb C^{|G|+1}$, for all $g\in G$, and $H^m:\cong\mathbb C^{2|G|}$. Let us fix a
basis $\{v^x_y\}(y\in\Omega_x)$ in each $H^x\ (x\in\Omega)$ choosing a norm one vector in every 1-dimensional vector subspace:
$v^g_h\in Hom(h,(h-g)\otimes g)$, $v^g_m\in Hom(m,m\otimes g)$, $v^m_g\in Hom(m,g\otimes m)$, and $v^m_{\overline g}\in
Hom(g,m\otimes m)$, where $g\in G$.
Now the whole WHA structure of $\mathfrak G_{\mathcal T\mathcal Y}=(B,\Delta,S,\varepsilon)$ is given by formulas
(\ref{algebra}), (\ref{coproduct}), (\ref{counit}),  (\ref{product}) and  (\ref{antipode}). In particular, the $C^*$-algebra
$B= \underset{x \in \Omega} \oplus H^x \otimes \overline{H^x}$ has a canonical basis  $\{f^x_{\alpha,\beta}=v^x_\alpha\otimes
\overline{v^x_\beta}\}_{x \in \Omega, \alpha,\beta \in \Omega_x}$.
\vskip 0.5cm
For all $x,y \in \Omega$ and all $v \in H^x$, $w \in H^y$, denote $ v\circ w = \mathcal J_{x,y}(v \otimes_R w)$.
Then for  all   $\alpha, \beta \in \Omega_x, \gamma, \delta \in \Omega_y$, one has:
$$
 f^x_{\alpha,\beta}f^y_{\gamma,\delta} = (v^x_ \alpha \circ
v^y_\gamma) \otimes \overline {(v^x_\beta \circ v^y_\delta}),
$$
where computations made in \cite{M} 2.1.5, give, for all   $g,h,k\in G$:
$$
v^g_k\circ v^h_{x}= \delta_{x,h+k}v^{g+h}_{h+k},\ v^g_m\circ v^h_x= \delta_{x,m}v^{g+h}_m,
$$
$$
v^m_k\circ v^g_x=\delta_{x,m}\chi(g,k)v^m_k,\ v^m_{\overline k}\circ v^g_{x}=
\delta_{x,g+k}v^m_{\overline{g+k}},
$$
$$
v^g_x\circ v^m_{\overline k}=\delta_{x,m}\chi(g,k) v^m_{\overline k},\ v^g_x\circ v^m_k
= \delta_{x,k}v^m_{k-g},
$$
$$
v^m_h\circ v^m_{\overline k}= v^{k-h}_k,
\ v^m_{\overline h}\circ v^m_k=\delta_{h,k}\tau \underset{g\in G}\Sigma \chi(g,h)^{-1}v^g_m.
$$
The coproduct and the counit are defined, respectively, by
$$
\Delta(f^x_{\alpha,\beta}) = \underset{\alpha',\beta' \in \Omega_x} \sum f^x_{\alpha,\alpha'}\otimes
f^x_{\beta',\beta}
$$
and $\varepsilon(f^x_{\alpha,\beta})=\delta_{\alpha,\beta}.$ The antipode and the involution are as follows:
%$\Phi_x:H^x\to\overline{H^{x^*}}$ and $\Psi_x:\overline{H^x}\to H^{x^*}$ are given, respectively, by:
%$$
%\Phi_g(v^g_h)=\overline{v^{-g}_{h-g}},\ \Phi_g(v^g_m)=\overline{v^{-g}_m},\ \Phi_m(v^m_g)=|G|^{-1/2}\overline{v^m_{\overline g}},\
%\Phi_m(v^m_{\overline g})=\tau|G|^{-1/2}\overline{v^m_g},
%$$
%and
%$$
%\Psi_g(\overline{v^g_h})=v^{-g}_{h-g},\ \Psi_g(\overline{v^g_m})=v^{-g}_m,\ \Psi_m(\overline{v^m_g})=|G|^{1/2}v^m_{\overline g},\
%\Psi_m(\overline{v^m_{\overline g}})=\tau^{-1}|G|^{1/2}v^m_g,
%$$
%and :
%$$
%(v^g_k)^\sharp = v^{-g}_{k-g}, (v^g_m)^\sharp = v^{-g}_m, (v^m_g)^\sharp = |G|^{1/2}v^m_{\overline g}, (v^m_{\overline g})^\sharp = \tau^{-1}|G|^{1/2}v^m_g
%$$
%which implies that:
\begin{equation} \label{pode1}
S(f^g_{h,k})=f^{-g}_{k-g,h-g},\ S(f^g_{h,m})=f^{-g}_{m,h-g},\ S(f^g_{m,h})=f^{-g}_{h-g,m},
\end{equation}
\begin{equation} \label{pode2}
S(f^g_{m,m})=f^{-g}_{m,m},\ S(f^m_{g,h})=f^m_{\overline h,\overline g},\ S(f^m_{g,\overline h})=\tau^{-1} f^m_{h,\overline g},
\end{equation}
$$
S(f^m_{\overline g,h})=\tau f^m_{\overline h,g},\ S(f^m_{\overline g,\overline h})=f^m_{h,g}
$$
and:
$$
(f^g_{h,k})^*=f^{-g}_{h-g,k-g},\ (f^g_{h,m})^*=f^{-g}_{h-g,m},\ (f^g_{m,h})^*=f^{-g}_{m,h-g},
$$
$$
(f^g_{m,m})^*=f^{-g}_{m,m},\ (f^m_{g,h})^*=f^m_{\overline g,\overline h},\ (f^m_{g,\overline h})^*=\tau f^m_{\overline g,h},
$$
$$
(f^m_{\overline g,h})^*=\tau^{-1} f^m_{g,\overline h},\ (f^m_{\overline g,\overline h})^*=f^m_{g,h}.
$$
Recall that $H^0$ is a commutative $C^*$-algebra isomorphic to $R \sim  \mathbb C^\Omega$.
\vskip 0.5cm
\begin{remark}
\label{camille}
Since $\mathfrak G_{\mathcal T\mathcal Y}$ is selfdual, we also have $B=\underset{g\in G}\oplus B^g\oplus B^m$, where
$B^g\cong M_{|G|+1}(\mathbb C),\ \forall g\in G,\ B^m\cong M_{2|G|}(\mathbb C)$ (see \cite{M}, 2.1). Using the basis
$\{f^x_{y,z}\}$ and the matrix units $\{e^x_{y,z}\}$ of $B$ with respect to the basis $\{v^x_t\}$ of $H^x$, any irreducible
corepresentation $U^x\ (x \in \Omega)$ of $\mathfrak G_{\mathcal T\mathcal Y}$ can be written as
$$
U^x = \underset{ y,z\in \Omega_x} \sum e^x_{y,z} \otimes f^x_{y,z}.
$$

%In \cite{M} 2.1 is given an alternative description of the quantum groupoid $\mathfrak G_{\mathcal T\mathcal Y}=(B,\Delta,S,\varepsilon)$
%as we have $B=\underset{g\in G}\oplus B^g\oplus B^m$, where $B^g\cong M_{|G|+1}(\mathbb C),\ \forall g\in G,\ B^m\cong M_{2|G|}(\mathbb C)$
%for special  type isomorphisms  (\cite{M}, Proposition 2.1.11). One will consider the matrix basis of $B$, namely $(e^x_{y,z})$,  associated
%with  vector basis $(v^x_t)$ of each $H^x$, to this basis is attached  a co-matrix basis namely $(f^x_{y,z})$.

%ii)  Family $(U^x = \underset{ y,z  \in \Omega_x} \sum e^x_{y,z} \otimes f^x_{y,z})_{ x \in \Omega}$ (with notations of \cite{M} 2.1.)
%is   in $Ucorep \mathfrak G_{\mathcal T\mathcal Y}$ and gives an exhaustive set of irreducibles.

\end{remark}

\subsection{Quantum subgroupoids and quotient type coideals}

\begin{remark} \label{lattice} The lattice $Subgrp(G)$ of subgroups of $G$ with operations $\cap$ and $\vee=+$
can be extended to $\overline{Subgrp(G)}:=Subgrp(G)\sqcup\{\Omega\}$, where $\Omega=G\sqcup\{m\}$, by putting $L\cap
\Omega=\Omega\cap L=L$ and $L\vee\Omega=\Omega\vee L=\Omega\cap\Omega=\Omega\vee\Omega=\Omega$, for any subgroup
$L$ of $G$. Any rigid tensor $C^*$-subcategory of $\mathcal T\mathcal Y(G,\chi,\tau)$ is equivalent either to $Vec_L$ -
the category of finite dimensional $L$-graded vector spaces $(L<G)$ or to $\mathcal T\mathcal Y(G,\chi,\tau)$.
Let $\mathcal C^x\ (x\in\overline{Subgrp(G)})$ be a representative in such equivalence class of subcategories, in
particular, $\mathcal C^\Omega=\mathcal T\mathcal Y(G,\chi,\tau)$.
%are parameterized by the elements of $\overline{Subgrp(G)}$.
\end{remark}

In order to construct all quotient type coideals of $\mathfrak G_{\mathcal T\mathcal Y}$, first construct all its quantum
subgroupoids (up to isomorphism). Theorems \ref{genreconstruction} and \ref{mapreconstr} imply that any quantum subgroupoid of
$\mathfrak G_{\mathcal T\mathcal Y}$ is isomorphic to one of the quantum subgroupoids $(\mathfrak G^x,\pi^x)$ such that
$URep(\mathfrak G^x)\cong\mathcal C^x$, where $x\in\overline{Subgrp(G)}$. Define
%and that any quantum subgroupoid of $\mathfrak G_{\mathcal T\mathcal Y}$ is Morita equivalent to one of them. Define
$(\mathfrak G^\Omega,\pi^\Omega)=(\mathfrak G_{\mathcal T\mathcal Y},id)$ and, for any $L<G$, $(\mathfrak G^L,\pi^L)$ as follows:

\begin{lemma}\label{paysagecomplet1} If $e^x_{y,z}\ (x\in\Omega,y,z\in\Omega_x)$ are the matrix units of $B$ (see Remark
\ref{camille}), $l \in L, g,h \in G$, the collection $(B_L,\Delta_L,S_L,\varepsilon_L)$, where $B_L = \underset{l \in L}\oplus B(H^l)$,
$$
\Delta_L(e^l_{g,g'}) = \underset{ \underset {l_1+l_2  = l} {l_1,l_2 \in L}} \sum e^{l_1}_{g-l_2,g'-l_2} \otimes e^{l_2}_{g,g'}, \    \Delta_L(e^l_{g,m}) =   \underset{ \underset {l_1+l_2  = l} {l_1,l_2 \in L}} \sum e^{l_1}_{g-l_2,m} \otimes e^{l_2}_{g,m},
$$
$$
\Delta_L(e^l_{m,g'})=\underset{ \underset {l_1+l_2  = l} {l_1,l_2 \in L}} \sum e^{l_1}_{m,g'-l_2} \otimes e^{l_2}_{m,g'},
\   \Delta_L(e^l_{m,m}) =   \underset{ \underset {l_1+l_2  = l} {l_1,l_2 \in L}} \sum e^{l_1}_{m,m} \otimes e^{l_2}_{m,m},
$$
$$
S_L(e^l_{g,g'}) =   e^{-l}_{g'-l,g-l}, \ \quad    S_L(e^l_{g,m}) =  e^{-l}_{m,g-l},
$$
$$
S_L(e^l_{m,g'}) =   e^{-l}_{g'-l,m}, \ \quad    S_L(e^l_{m,m}) =  e^{-l}_{m,m},
$$
$$
\varepsilon_L(e^l_{g,g'}) = \varepsilon_L(e^l_{m,g'}) =   \varepsilon_L(e^l_{g,m'}) =\varepsilon_L(e^l_{m,m}) = \delta_{l,0},
$$
defines a WHA $\mathfrak G^L$. The canonical projection $\pi_L :B\to B_L$ defined, for all $x \in\Omega,\alpha,\beta\in\Omega_x$,
by $\pi_L(e^x_{\alpha,\beta}) = \delta_{x,L}e^x_{\alpha,\beta}$, where $\delta_{x,L} = 1$ if $x\in L$ and $=0$ otherwise, gives
to ${\mathfrak G}^L$ the structure of a quantum subgroupoid of $\mathfrak G_{TY}$.
\end{lemma}
\begin{dm}
Straightforward computations.
\hfill\end{dm}
%then $\mathfrak G^L = (\underset{l \in L} \oplus End(H^l),  \Delta_L,  S_L,  \varepsilon_L)$ is a $C^*$-quantum groupoid ($*$-WHA).
\begin{corollary} \label{camille1}
A linear basis for $(B_L)_t$ (resp., $(B_L)_s$) is given by $(e(L)^\alpha)_{\alpha \in \Omega}$ (resp.,
$(e(L)_\alpha)_{\alpha \in \Omega}$), where, for all $g \in G$, one has: $ e(L)^g =\underset{l \in L}\sum e^l_{g+l,g+l}$,
$e(L)_g = \underset{l \in L} \sum e^l_{g,g}$ and $ e(L)^m =  e(L)_m = \underset{l \in L} \sum e^l_{m,m}$.

The counital maps are given by:
$$
\epsilon_t^{B_L}(e^l_{x,y}) = \delta_{l,0} e(L)^x,\quad  \epsilon_s^{B_L}(e^l_{x,y}) = \delta_{l,0} e(L)_y.
$$
A linear basis for $(B_L)_s\cap (B_L)_t$ is given by $(z_\beta)_{\beta \in G/L \sqcup \{m\}}$, where, for all
$\beta \in G/L$, one has: $z_\beta = \underset{l \in L,g,\in \beta} \sum e^l_{g,g}$ and $z_m =  e(L)^m =  e(L)_m =  \underset{l
\in L} \sum e^l_{m,m}$.
Moreover, $(B_L)_s\cap Z(B_L) = \mathbb C $, so ${\mathfrak G}^L$ is connected and not coconnected.
\end{corollary}

\begin{remark} \label{hopf} Any ${\mathfrak G}^L$ is Morita equivalent to a commutative and cocommutative Hopf
$C^*$-algebra generated by the group $L$.
\end{remark}

\begin{proposition}
\label{puma1}
%For all $x \in \overline {Subgrp(G)}$, $ \mathfrak G^x$ is a quantum subgroupoid of $\mathfrak G_{\mathcal T\mathcal Y}$.
Denote $I^x := I(\mathfrak G^x\backslash\mathfrak G)$. Then $I^\Omega = B_s$ and, for any subgroup $L$ of $G$, setting
$v^0_Y:= \underset{y \in Y} \sum v^0_y$, where $Y\in G/L$, one has:
$$
I^L =   Vec <v^0_Y, Y \in G/L> \otimes \overline {H^0} \underset{l \in ÊL^\perp} \oplus  v^l_m \otimes \overline {H^l}.
$$
%$$I^\Omega = B_s$$
\end{proposition}
\begin{dm}  We will use Lemma \ref{genregular}.
%For all subgroup $L \subset G$, if  $\pi_L :B  \to  {\mathfrak G}^L$ is the canonical projection given for all $ x \in
%\Omega, \alpha ,\beta \in \Omega_x$ by $\pi_L(e^x_{\alpha,\beta}) = \delta_{x,L}e^x_{\alpha,\beta}$ where $\delta_{x,L} =
%1$ when $x \in L$ and $=0$ otherwise.
%Then  map $\pi_L$  gives to  ${\mathfrak G}^L$ the structure of  a quantum subgroupoid of $\mathfrak G_{TY}$.
For all $\hat c \in \hat {\mathfrak G}^L, \beta \in \Omega_\alpha$,  one has:
\begin{equation}
\label{soucis}
(\pi_L)_*(\hat c)v^\alpha_\beta
= \underset{i,j \in \Omega_\alpha} \sum <f^\alpha_{i,j}, (\pi_L)_*(\hat c)>e^\alpha_{i,j}v^\alpha_\beta
= \underset{i\in \Omega_\alpha} \sum <\pi_L(f^\alpha_{i,\beta}),  \hat c> v^\alpha_i
\end{equation}
A linear form $\phi$  on $B_L$ is a left integral if and only if $(i\otimes \phi)\Delta_L(e^l_{x,y})$ is in $(B_L)_t$, for all
$l\in L,x,y \in \Omega$. Then Lemma \ref{paysagecomplet1} implies that $\Lambda =  Vec<\lambda_x, x \in \Omega>$.

For all $\alpha,\gamma \in \Omega,  \beta \in \Omega_\alpha$, using (\ref{soucis}), one has:
\begin{align*}
(\pi_L)_*(\lambda_\gamma)v^\alpha_\beta= \underset{i\in \Omega_\alpha, l \in L} \sum <\pi_L(f^\alpha_{i,\beta}), (e^l_{\gamma,\gamma})^*>
v^\alpha_i.
\end{align*}
This gives, in particular:
%For $\alpha = p \in G$:
%\begin{align*}
%(\pi_L)_*(\lambda_\gamma)v^p_\beta= \underset{h \in G, l \in L} \sum <\pi_L(f^p_{h,\beta}), (e^l_{\gamma,\gamma})^*>
%v^p_h+ <\pi_L(f^p_{m,\beta}), (e^l_{\gamma,\gamma})^*> v^p_m
%\end{align*}
%In particular:
%- if $\beta = k \in G:  \ \ \ \ \
$$
(\pi_L)_*(\lambda_\gamma)v^p_k= \delta_{p,0} \delta_{ \gamma,-k}\underset{l \in L} \sum v^0_{k+l},
$$
%- if $\beta = m: \ \ \ \ \ \ \ \ \ \
$$
(\pi_L)_*(\lambda_\gamma)v^p_m= \delta_{\gamma,m} \delta_{p,L^\perp} |L|v^p_m
$$
And also:
%\vskip 0.5cm
%For $\alpha = m$:
%\begin{align*}
%(\pi_L)_*(\lambda_\gamma)v^m_\beta= \underset{h \in G, l \in L} \sum <\pi_L(f^m_{h,\beta}), (e^l_{\gamma,\gamma})^*> v^m_h+\underset{h \in G, l \in L} \sum <\pi_L(f^m_{\overline h,\beta}), (e^l_{\gamma,\gamma})^*> v^m_{\overline h}
%\end{align*}
%In particular:
%- if $\beta = k \in G:
$$
(\pi_L)_*(\lambda_\gamma)v^m_k= \underset{h \in G, l \in L} \sum <\pi_L(e^{h-k}_{-k,m}), (e^l_{\gamma,\gamma})^*> v^p_h= 0,
$$
%- if $\beta = \overline k \in \overline G:
$$
(\pi_L)_*(\lambda_\gamma)v^m_{\overline k}
= \underset{h \in G, l \in L} \sum  <\pi_L(e^{h-k}_{m,-k}), (e^l_{\gamma,\gamma})^*> v^m_{\overline h}  =  0.
$$
\vskip 0.5cm
So if one sets: $v^0_X = \underset{x \in X} \sum v^0_x$, for all $X \subset \Omega$, then:
$$
(\pi_L)_*(\Lambda) H^p = \delta_{p,0}( Vec<v^0_Y, v^0_m / Y \in G/L>)+ \delta_{p, L^\perp} \mathbb C v^p_m\ \text{for\ all}\ p \in G,
$$
$$
(\pi_L)_*(\Lambda) H^m = \{0\}.
$$
These calculations and Lemma \ref{genregular} give the result.
\hfill \end{dm}
\end{subsection}

\begin{subsection} {The lattice of invariant  coideals}

In order to precise the relationship between quotient type and invariant coideals and to
characterize the lattice of these coideals in the Tambara-Yamagami case, rewrite the definition
(\ref{adjoint}) of $\triangleleft$ using (\ref{pode1}) and (\ref{pode2}) as follows:

\begin{align*}
(\eta^y\otimes\overline {\xi^y})\triangleleft(\eta^x\otimes\overline{\xi^x})
=(\underset{z\in\Omega_x}\sum(v^x_z)^\natural\circ\eta^y\circ v^x_z)\otimes(\eta^x)^\flat\circ
\overline{\xi^y }\circ\overline{\xi^x},
\end{align*}
where $x,y \in \Omega, \eta^x,\xi^x \in H^x, \eta^y,\xi^y \in H^y$. This expression allows to define
the map $P^x:\underset{y\in\Omega}\oplus H^y\to\underset{y\in\Omega}\oplus H^{y}$ by putting
for any fixed $x,y\in\Omega,\eta^y\in H^x$:
$$
P^x(\eta^y)=\underset{z\in\Omega_x}\sum(v^x_z)^\natural\circ \eta^y\circ v^x_z,
$$
%does not depend on the choice of $\{e^x_i\}$, so it
and we have:
\begin {lemma}\label{caniculerecord}
A  coideal $I=\underset{y\in\Omega}\sum X^y \otimes\overline{H^y}$ is invariant if and only if
$\underset{y\in\Omega}\oplus X^y =P^x(\underset{y\in\Omega}\oplus X^y)$ for all $x\in\Omega$.
\end{lemma}
A straightforward calculation of $P^x$ on the basic elements $v^x_z$ proves
\begin{lemma} \label{neglige}
For all $g,k,h\in G$, one has:
\begin{align*}
P^h( v^g_k ) = \delta_{g,0}v^0_{k+h}
 \
 &, \
P^m(v^g_k) = \delta_{g,0}sign(\tau) \underset{p \in G} \sum \chi(p,k)v^p_m  \\
P^h( v^g_m ) = v^g_m
 \
 &, \
P^m(v^g_m ) = \tau^{-1}|G|^{1/2}\underset{p \in  G} \sum  \chi(g,p)  v^0_p \\
P^h_{\mid{H^m}}
&=
P^m_{\mid{H^m}}
  =0
 \end{align*}
 \end{lemma}

In the Tambara-Yamagami case any invariant coideal is not only isomorphic but is itself of quotient type:
\begin{proposition}
\label{puma2}
For any invariant coideal $I= \underset{y \in \Omega} \oplus X^y \otimes \overline {H^y}$ there is a unique quotient type
coideal $I^x = I(\mathfrak G^x  \backslash\mathfrak G)$ such that $I=I^x$.

%The map: $x \mapsto I^x$ is a decreasing anti-isomorphism of lattices between $\overline{Subgr(G)}$ and $Invl(B)$.
\end{proposition}
\begin{dm}
%First we prove that any invariant coideal $I= \underset{y \in \Omega} \oplus X^y \otimes \overline {H^y}$ equals
%and not just isomorphic to $I^x = I(\mathfrak G^x  \backslash\mathfrak G)$.
Due to \cite{VV2}, Lemma 3.3, b) there is a partition $(\Gamma_i)_{i \in I^0}$ of $\Omega$ such that $X^0 = Vec<v^0_{\Gamma_i}
/ i \in I^0>$. Moreover, putting $K :=\{g \in G / Dim(X^g) \not= 0\}$, one has by Lemmas \ref{caniculerecord} and \ref{neglige}:
$X^m = \{0\}$, and $X^g = \mathbb Cv^g_m$ for all $g \in K$, $g\not= 0$.
With the convention $m+g=m$ for all $g \in G:$, one has by Lemma \ref{neglige}: $v^0_{\Gamma_i+g} \in X^0$ for all $i \in I^0$.
As a consequence, for all $g \in G$ there is $j \in J$ such that $\Gamma_i+g = \Gamma_j$; this allows only two possibilities:

1) there is a single class $\Gamma_i = \Omega$, so $I = B_s+ \underset{g \in K\setminus \{0\} }\oplus \mathbb Cv^g_m$. But for all
$g \in K\setminus \{0\}$ one has $P^m(v^g_m ) = \tau^{-1}|G|^{1/2}\underset{p\in G} \sum\chi(g,p)v^0_p$ which must be collinear to
$v^0_\Omega=\underset {y\in\Omega}\sum v^0_y$. Hence, $K= \{0\}$ and $I= B_s=I^\Omega$.

2) The partition $(\Gamma_i)$ is $\{m\}$ together with a partition $(\Gamma^G_p)_{p \in P}$ of $G$. Moreover, for any $p,q \in P$
there is $g \in G$ such that $\Gamma^G_p+g = \Gamma^G_q$. For any $p \in P$ denote $L_p = \{g\in G/\Gamma^G_p+g =\Gamma^G_p\}$,
then $L_p$ is a subgroup of $G$. But since for all $q \in P$, there is $h \in G$ such that $\Gamma_q = \Gamma_p+h$, the group
$L_p$ does not depend on $p \in P$, denote it by $L$. Let us show that $\Gamma^G_p \in G/L$ for any $p\in P$.

If $h\in\Gamma^G_p - \Gamma^G_p$, then $(\Gamma^G_p + h)\cap\Gamma^G_p\not=
\emptyset$, so $\Gamma^G_p + h  = \Gamma^G_p $. Hence, $K =\Gamma^G_p - \Gamma^G_p$.
For all $p \in P$, let $z$ be in $\Gamma^G_p$, obviously one has $z+L\subset \Gamma^G_p$, let $t \in \Gamma^G_p$, then $t-z \in L=
\Gamma^G_p - \Gamma^G_p$ hence $t \in z+L$, as a consequence $z+L = \Gamma^G_p$ so $\Gamma^G_p \in G/L$, as we deal with a partition of
$G$: $\{ \Gamma^G_p/ p \in P\} = G/L$. So we have: $X^0 = Vec<(v^0_Y)_{Y \in G/L},v^0_m>$.

If now $g \in K$, then due to Lemma \ref{neglige} one must have $\underset{h \in G} \sum\chi(g,h)v^{0}_{h } \in X^0$, but:
$ \underset{h \in G} \sum \chi(g,h)v^{0}_{h } = \underset{p\in G/L} \sum \underset{h \in p} \sum \chi(g,h)v^{0}_{h }$,
so this element has to belong to $X^0$ and must be of the form $\underset{p\in G/K}\sum  \mu_p  v^0_{p}$, i.e., for all $p \in G/K$
and all $h \in p$, one has $\mu_p = \chi(g,h)$, which means that $K\subset L^\perp$. Conversely, by Lemma \ref{neglige}, one must
have $\underset{k \in K, p \in G} \sum \chi(p,k)v^p_m \in \underset{k \in G} \oplus A^k$, but on the other hand:

\begin{align*}
\underset{k \in K} \sum  \underset{p \in G} \sum \chi(p,k)v^p_m  =  \underset{p \in G} \sum (\underset {k \in K} \sum \chi(k,p))v^p_m = |K| \underset{p \in K^\perp } \sum v^p_m
\end{align*}
So  $L^\perp \subset K$ and in case 2) we have $I= I^L$.
%Let us prove that the map: $x \mapsto I^x$ is a bijection between the sets $\overline{Subgr(G)}$ and $Invl(B)$.
%By Lemma \ref{quinv},  $I^x  = I(\mathfrak G^x  \backslash\mathfrak G)$ is invariant, so we have to prove that all invariant
%coideal $I= \underset{y \in \Omega} \oplus X^y \otimes \overline {H^y}$ is of the form $I^x$.
\hfill\end{dm}
\begin{corollary} \label{big}
In $\mathfrak G_{\mathcal T\mathcal Y}(G,\chi,\tau)$, the sets of quotient type and invariant coideals   coincide and
are in bijection with $\overline{Subgrp(G)}$.
%\end{corollary}%The  map: $x \mapsto I^x$ is a bijection $\overline{Subgrp(G)}\to Invl(B)$.
\end{corollary}
\begin{dm} By \ref{puma2} any invariant coideal is quotient type, conversely any quotient type coideal is invariant by lemma \ref{quinv}, moreover the set  $\{I^x, x \in \overline{Subgrp(G)}\}$ contains all invariant coideals and is included in the set of quotient type coideals.
\hfill \end{dm}

\vskip 0.5cm 

Finally, we describe the lattice of invariant (or quotient type) coideals.
\begin{proposition}
\label{puma3}
The map: $x \mapsto I^x$ is an anti-isomorphism of the lattices $\overline{Subgr(G)}$ and $Invl(B)$.
\end{proposition}
\begin{dm}
For any subgroup $K$ of $G$, since $I^K = Vec <v^0_Y, Y \in G/K> \otimes \overline {H^0} \underset{k \in ÊK^\perp} \oplus v^k_m
\otimes \overline {H^k}$, one sees that the map $x \mapsto I^x$ is decreasing. Hence, for all $L,K < G$, one has $ I^{L +K}\subset
I^L \cap I^K$. Conversely, for all $z \in I^L \cap I^K $, there exist some families of complex numbers $(\lambda_u),(\mu_v)$ and
non zero vectors  $(\xi^0_w), (\eta^g_m)$ such that:
$$
z = \underset{Y  \in G/L} \sum \lambda_Y v^0_Y \otimes \overline {\xi^0_Y }+ \underset{ l \in L^\perp} \sum \lambda_l v^l_m
\otimes \overline {\eta^l_m}=  \underset{Z  \in G/K} \sum \mu_Z v^0_Z \otimes \overline {\xi^0_Z }+ \underset{ k \in K^\perp}
\sum \mu_k v^k_m \otimes \overline {\eta^k_m}.
$$
This gives:
$$
\underset{g \in G} \sum \lambda_{\overset .g} v^0_g \otimes \overline {\xi^0_{\overset .g} }+ \underset{ l \in L^\perp}
\sum \lambda_l v^l_m \otimes \overline {\eta^l_m}=  \underset{g \in G} \sum \lambda_{\overset {..}g} v^0_g \otimes \overline
{\xi^0_{\overset {..}g} }+ \underset{ k \in K^\perp} \sum \lambda_k v^k_m \otimes \overline {\eta^k_m},
$$
where $\overset . g$ (resp., $\overset{..}g$) is the class of $g$ in $G/L$ (resp., in $G/K$). As a consequence, one has:

1) if $l \in L^\perp$ and $l\notin K^\perp$, then  $\lambda_l = 0$.

2) for any $g \in G$ the following equality holds:  $\lambda_{\overset .g} \overline {\xi^0_{\overset .g} }=\lambda_{\overset {..}g} \overline
{\xi^0_{\overset {..}g} }.$

Condition 2) implies that for all $g \in G, p \in L +K$, one has $\lambda_{\overset .g} \overline {\xi^0_{\overset .g} } =
\lambda_{\overset .{\overline{ g+p}}} \overline {\xi^0_{\overset .{\overline{ g+p}}} }$. Hence, for any $Y \in G/(L+ K)$ we
can define $\overline{\xi^0_Y}$ such that for all $g \in Y$ one has: $\lambda_{\overset .g} \overline {\xi^0_{\overset .g} } =
\overline{\xi^0_Y}$. Then, using the fact that $L^\perp\cap K^\perp= (L+K)^\perp$ (see \cite{HR}, (23),(29)(b) on p.369), one has:
\begin{align*}
z
&= \underset{g \in G} \sum \lambda_{\overset .g} v^0_g \otimes \overline {\xi^0_{\overset .g} }+ \underset{ l \in L^\perp
\cap K^\perp} \sum \lambda_l v^l_m \otimes \overline {\eta^l_m} \\
&= \underset{Y \in  G/(L+K)} \sum \  \underset{g\in Y} \sum   v^0_g \otimes \overline {\xi^0_{Y} }+ \underset{ l \in (L+K)^\perp}
\sum \lambda_l v^l_m \otimes \overline {\eta^l_m}\\
&= \underset{Y \in  G/(L+K)} \sum   v^0_Y \otimes \overline {\xi^0_{Y} }+ \underset{ l \in (L+K)^\perp} \sum \lambda_l v^l_m
\otimes \overline {\eta^l_m}
\end{align*}
Hence, $z \in I^{L+ K}$, which proves that $ I^{L +K} = I^L \cap I^K $.

Obviously, $I^L\vee I^K\subset I^{L \cap K}$. Conversely, since $(L\cap K)^\perp = L^\perp+ K^\perp$, for any
$p\in (L\cap K)^\perp$ there exist $l\in L^\perp$ and $k\in K^\perp$ such that $p = k+l$. Then $v^{p}_m \otimes
\overline{ H^{p}} =(v^{l}_m\otimes\overline{ H^{l}})(v^{k}_m \otimes \overline{ H^{k}})$, so it belongs to
$I^L\vee I^K$. For all $g \in G$ one has:
$$
v^0_{g+L} \circ v^0_{g+K}= v^0_{g+L\cap K},
$$
hence, $v^0_{g+L\cap K}\otimes\overline{H^0}$ belongs to $I^L\vee I^K$. All basic elements of $I^{L\cap K}$ are in
$I^L\vee I^K$ which gives the converse inclusion $ I^{L\cap K}\subset I^L\vee I^K$, and the result follows.
\hfill \end{dm}
\end{subsection}
\end{section}

%By Lemma \ref{quinv} $I^x$ is invariant.

%XXXXXXXXXXXXXXXXXXXXXXXXXXXXXXXXXXXXXXXXXXXXXXXXXXXXXXXXXXXXXXXXXXXXXXXXXXXXXx

%So we know, using Proposition \ref{puma1}, Lemma \ref{quinv} and Remark \ref{coronalist}, that $ \{I^x / x  \in \overline{Subgrp(G)}\}$
%is an exhaustive subset of $Invl(B)$. Let us prove that, more precisely, $Invl(B) = \{I^x / x  \in \overline{Subgrp(G)}\}$.

%\begin{corollary}
%\label{repos}
%Map: $x \mapsto I^x$ is a decreasing anti-isomorphism of lattices between $\overline{Subgr(G)}$ and $Invl(B)$.
%\end{corollary}

%\end{subsection}
%\end{section}

%%%%%%%%%%%%%%%%%%%%%%%%%%%%%%%%%%%%%%%%%%%%%%%%%%%%%%%%%%%%%%%%%%%%%%%%%%%%%%%%%%%%%%%%%%%%%%%

%%%%%%%%%%%%%%%%%%%%%%%%%%%%%%%%%%%%%%%%%%%%%%%%%%%%%%%%%%%%%%%%%%%%%%%%%%%%%%%%%%%%

\bibliographystyle{plain}
\bibliography{biblio1}
\end{document}